\documentclass[11pt]{article}
\usepackage{amsmath}
\usepackage{color}
\newcommand{\Limsup}{\mathop{{\rm Lim}\,{\rm sup}}}
\def\disp{\displaystyle}

\def\tto{\;{\lower 1pt \hbox{$\rightarrow$}}\kern -10pt
\hbox{\raise 2pt \hbox{$\rightarrow$}}\;}
\def\Hat{\widehat}
\def\hat{\widehat}
\def\Tilde{\widetilde}
\def\tilde{\widetilde}
\def\Bar{\overline}
\def\ra{\rangle}
\def\la{\langle}
\def\ve{\varepsilon}
\def\B{I\!\!B}
\def\IN{I\!\!N}
\def\N{\mathcal{N}}

\def\h{\hfill\Box}
\def\R{I\!\!R}

\def\ox{\bar{x}}

\def\co{\mbox{\rm co}\,}

\def\epi{\mbox{\rm epi}\,}
\def\dim{\mbox{\rm dim}\,}
\def\dom{\mbox{\rm dom}\,}

\def\cl{\mbox{\rm cl}^*}
\def\cone{\mbox{\rm cone}\,}

\def\substack#1#2{{\scriptstyle{#1}\atop\scriptstyle{#2}}}

\def\h{\hfill\triangle}
\def\dn{\downarrow}
\def\O{\Omega}

\def\ph{\varphi}
\def\emp{\emptyset}
\def\st{\stackrel}
\def\oR{\Bar{\R}}
\def\tR{\Tilde\R^T_+}
\def\lm{\lambda}
\def\gg{\gamma}
\def\dd{\delta}

\newcounter{lk}
\def\Limsup{\mathop{{\rm Lim}\,{\rm sup}}}

\begin{document}

\newtheorem{Theorem}{Theorem}[section]
\newtheorem{Proposition}[Theorem]{Proposition}
\newtheorem{Remark}[Theorem]{Remark}
\newtheorem{Lemma}[Theorem]{Lemma}
\newtheorem{Corollary}[Theorem]{Corollary}
\newtheorem{Definition}[Theorem]{Definition}
\newtheorem{Example}[Theorem]{Example}
\renewcommand{\theequation}{\thesection.\arabic{equation}}
\normalsize

\def\proof{
\normalfont
\medskip
{\noindent\itshape Proof.\hspace*{6pt}\ignorespaces}}
\def\endproof{$\h$ \vspace*{0.1in}}
\begin{center}
\vspace*{0.3in} {\bf CONSTRAINT QUALIFICATIONS AND OPTIMALITY CONDITIONS FOR NONCONVEX SEMI-INFINITE AND INFINITE
PROGRAMS}\footnote{Research was partially supported by the USA National Science Foundation under grants
DMS-0603846 and DMS-1007132 and by the Australian Research Council under grant DP-12092508.}\\[3ex]
BORIS S. MORDUKHOVICH \footnote{Department of Mathematics, Wayne
State University, Detroit, Michigan 48202; email:
boris@math.wayne.edu.} and T. T. A. NGHIA\footnote{Department of
Mathematics, Wayne State University, Detroit,
Michigan 48202; email: ttannghia@gmail.com.}\\[2ex]
{\bf Dedicated to Jon Borwein in honor of his 60th  birthday}\\[2ex]
\end{center}
\small {\bf Abstract.} The paper concerns the study of new classes
of nonlinear and nonconvex optimization problems of the so-called
infinite programming that are generally defined on
infinite-dimensional spaces of decision variables and contain
infinitely many of equality and inequality constraints with
arbitrary (may not be compact) index sets. These problems reduce
to semi-infinite programs in the case of finite-dimensional spaces
of decision variables. We extend the classical
Mangasarian-Fromovitz and Farkas-Minkowski constraint
qualifications to such infinite and semi-infinite programs. The
new qualification conditions are used for efficient computing the
appropriate normal cones to sets of feasible solutions for these
programs by employing advanced tools of variational analysis and
generalized differentiation. In the further development we derive
first-order necessary optimality conditions for infinite and
semi-infinite programs, which are new in both finite-dimensional
and infinite-dimensional settings.

\normalsize

\section{Introduction}
\setcounter{equation}{0}

The paper mainly deals with constrained optimization problems formulated as follows:
\begin{eqnarray}\label{1.1}
\left\{\begin{array}{ll}
\mbox{minimize }\;f(x)\;\mbox{ subject to}\\
g_t(x)\leq 0\;\mbox{ with }\;\;t\in T\;\mbox{ and }\;h(x)=0,
\end{array}\right.
\end{eqnarray}
where $f:X\to \Bar\R:=(-\infty,\infty]$ and $g_t:X\to \Bar \R$ as
$t\in T$ are extended-real-valued functions defined on Banach
space $X$, and where $h:X\to Y$ is a mapping between Banach
spaces. An important feature of problem (\ref{1.1}) is that the
index set $T$ is {\em arbitrary}, i.e., may be infinite and also
noncompact. When the spaces $X$ and $Y$ are finite-dimensional,
the constraint system in (\ref{1.1}) can be formed by finitely
many equalities and infinite inequalities. These optimization
problems belong to the well-recognized area of {\em semi-infinite
programming} (SIP); see, e.g., the books \cite{GL,GL1} and the
references therein. When the dimension of the decision space $X$
as well as the cardinality of $T$ are infinite, problem
\eqref{1.1} belongs to the so-called {\em infinite programming};
cf.\ the terminology in \cite{an,DGLS} for linear and convex
problems of this type. We also refer the reader to more recent
developments \cite{CLMP1,CLMP2,DMN1,DMN2,FLN,LNP} concerning
linear and convex problems of infinite programming with inequality
constraints. \vspace*{0.05in}

To the best of our knowledge, this paper is the first one in the
literature to address {\em nonlinear} and {\em nonconvex} problems
of infinite programming. Our primary goal in what follows is to
find verifiable {\em constraint qualifications} that allow us to
establish efficient {\em necessary optimality conditions} for
local optimal solutions to nonconvex infinite programs of type
\eqref{1.1} under certain differentiability assumptions on the
constraint (while not on the cost) functions. In this way we
obtain a number of results, which are new not only for infinite
programs, but also for SIP problems with noncompact (e.g.,
countable) index sets.

It has been well recognized in semi-infinite programming that the
{\em Extended Mangasarian-Fromovitz Constraint Qualification}
(EMFCQ), first introduced in \cite{JTW}, is particularly useful
when the index set $T$ is a compact subset of a finite-dimensional
space and when $g(x,t):=g_t(x)\in\mathcal{C}(T)$ for each $x\in
X$; see, e.g., \cite{BS,CHY,jrs,HK,KH,LS,S,st,zy} for various
applications of the EMFCQ in semi-infinite programming. Without
the compactness of the index set $T$ and the continuity of the
inequality constraint function $g(x,t)$ with respect to the index
variable $t$, problem (\ref{1.1}) changes dramatically and--as
shown below--does not allow us to employ the EMFCQ condition
anymore. That motivates us to seek for new qualification
conditions, which are more appropriate in applications to infinite
programs as well as to SIP problems with noncompact index sets and
infinite collections of inequality constraints defined by
discontinuous functions. \vspace*{0.05in}

In this paper we introduce two new qualification conditions, which
allow us to deal with infinite and semi-infinite programs of type
\eqref{1.1} without the convexity/linearity and compactness
assumptions discussed above. The first condition, called the {\em
Perturbed Mangasarian-Fromovitz Constraint Qualification} (PMFCQ),
turns out to be an appropriate counterpart of the EMFCQ condition
for infinite and semi-infinite programs \eqref{1.1} with
noncompact index sets $T$ and discontinuous functions
$g(x,\cdot)$. The second condition, called the {\em Nonlinear
Farkas-Minkowski Constraint Qualification} (NFMCQ), is a new
qualification condition of the {\em closedness} type, which is
generally independent of both EMFCQ and PMFCQ conditions even for
countable inequality constraints in finite dimensions.

Our approach is based on advanced tools of variational analysis
and generalized differentiation that can be found in \cite{M1,M2}.
Considerably new ingredients of this approach relate to computing
appropriate {\em normal cones} to the set of {\em feasible
solutions} for the infinite/semi-infinite program \eqref{1.1}
given by
\begin{equation}\label{1.2}
\O:=\big\{x\in X\big|\;h(x)=0,\;g_t(x)\le 0\;\mbox{ as }\;t\in
T\big\}.
\end{equation}

Since the feasible solution set $\O$ is generally nonconvex, we
need to use some normal cone constructions for nonconvex sets. In
this paper we focus on the so-called {\em Fr\'echet/regular normal
cone} and the {\em basic/limiting normal cone} introduced by
Mordukhovich; see \cite{M1} with the references and commentaries
therein. Developing general principles of variational analysis, we
employ this approach to  derive several necessary optimality
conditions for the class of nonlinear infinite programs under
consideration.\vspace*{0.05in}

The rest of the paper is organized as follows. In Section~2 we
present basic definitions as well as some preliminaries from
variational analysis and generalized differentiation widely used
in this paper. Section~3 is mainly devoted to the study of the new
PMFCQ and NFMCQ conditions for infinite programs in Banach spaces.
Relationships between the new qualification conditions and other
well-recognized constraint qualifications for SIP and infinite
programs are discussed here.

In Section~4, we provide exact computations for the Fr\'echet and
limiting normal cones to the feasible set of \eqref{1.1} under the
PMFCQ and NFMCQ conditions. This part plays a crucial role for the
subsequent results of the paper. Following this way, Section~5 concerns the
derivation of necessary optimality conditions for local minimizers
of the infinite and semi-infinite programs under consideration.
\vspace*{0.05in}

Our notation and terminology are basically standard and
conventional in the area of variational analysis and generalized
differentials.; see, e.g., \cite{M1,rw}. As usual, $\|\cdot\|$
stands for the norm of Banach space $X$ and $\la\cdot,\cdot\ra$
signifies for the canonical pairing between $X$ and its
topological dual $X^*$ with the symbol $\st{w^*}\to$ indicating
the convergence in the weak$^*$ topology of $X^*$ and the symbol
$\cl$ standing for the weak$^*$ topological closure of a set. For
any $x\in X$ and $r>0$, denote by $\B_r(x)$ the closed ball
centered at $x$ with radius $r$ while $\B_X$ stands for the closed
unit ball in $X$.

Given a set $\O\subset X$, the notation  $\co\O$ signifies the
convex hull of $\O$ while that of $\cone\O$ stands for the {\em
convex conic hull} of $\O$, i.e., for the convex cone generated by
$\O\cup\{0\}$.  Depending on the context, the symbols
$x\st{\O}{\to}\ox$ and $x\st{\ph}{\to}\ox$ mean that $x\to\ox$
with $x\in\O$ and $x\to\ox$ with $\ph(x)\to\ph(\ox)$ respectively.
Given finally a set-valued mapping $F\colon X\tto X^*$ between $X$
and $X^*$, recall that the symbol
\begin{eqnarray}\label{1.4}
\Limsup_{x\rightarrow \ox} F(x):=\Big\{x^*\in
X^*\Big|\;\exists\,x_n\to\ox,\;\exists\,
x^*_n\xrightarrow{w^*}x^*\;\mbox{ with }\;x^*_n\in F(x_n),\quad
n\in\IN\Big\}
\end{eqnarray}
stands for the {\em sequential Painlev\'{e}-Kuratowski outer/upper
limit} of $F$ as $x\to\ox$ with respect to the norm topology of
$X$ and the weak$^*$ topology of $X^*$, where
$\IN:=\{1,2,\ldots\}$.

\section{Preliminaries from Generalized Differentiation}
\setcounter{equation}{0}

In this preliminary section we briefly review some constructions
of generalized differentiation used in what follows; see
\cite{BZ,M1,rw,sc} for more details and related material.
Throughout this paper, unless otherwise stated, all the spaces
under consideration are Banach.

Given an extended-real-valued function $\ph\colon
X\to\oR:=(-\infty,\infty]$, we always assume that it is proper,
i.e., $\ph\not\equiv\infty$. The notation
\begin{eqnarray*}
\dom\ph:=\big\{x\in X\big|\;\ph(x)<\infty\big\}\;\mbox{ and }\;\epi\ph:=
\big\{(x,r)\in X\times\R\big|\;r\ge\ph(x)\big\}
\end{eqnarray*}
are used for the domain and the epigraph of $\ph$, respectively,

Define the {\em analytic $\ve$-subdifferential} of $\ph$ at
$\ox\in\dom\ph$ by
\begin{eqnarray}\label{2.1}
\Hat\partial_\ve\ph(\ox):=\Big\{x^*\in
X^*\Big|\;\liminf_{x\to\ox}\frac {\ph(x)-\ph(\ox)-\la
x^*,x-\ox\ra}{\|x-\ox\|}\ge-\ve\Big\},\quad\ve\ge 0
\end{eqnarray}
and let $\Hat\partial_\ve\ph(\ox):=\emp$ for $\ox\notin\dom\ph$.
If $\ve=0$, the construction $\Hat\partial
\ph(\ox):=\Hat\partial_0\ph(\ox)$ in (\ref{2.1}) is known as the
{\em Fr\'echet} or {\em regular subdifferential} of $\ph$ at
$\ox$; it reduces in the convex case to the classical
subdifferential of convex analysis. The sequential regularization
of \eqref{2.1} defined via the outer limit \eqref{1.4} by
\begin{eqnarray}\label{2.2}
\partial\ph(\ox):=\Limsup_\substack{x\st{\ph}{\to}\ox}{\ve\dn 0}\Hat\partial_\ve\ph(x),
\end{eqnarray}
is known as the {\em limiting}, or {\em basic}, or {\em
Mordukhovich subdifferential} of $\ph$ at $\ox\in\dom\ph$. It can
be equivalently described with $\ve=0$ in (\ref{2.2}) if $\ph$ is
lower semicontinuous (l.s.c.) around $\ox$ and if $X$ is an
Asplund space, i.e., each of its separable subspace has a
separable dual (in particular, any reflexive space is Asplund; see, e.g., \cite{BZ, M1}
for more details and references). We have $\partial\ph(\ox)\ne\emp$ for every locally Lipschitzian
function on an Asplund space.

A complementary construction to \eqref{2.2}, known as the {\em
singular} or {\em horizontal subdifferential} of $\ph$ at $\ox$,
is defined by
\begin{eqnarray}\label{2.3}
\partial^\infty\ph(\ox):=\Limsup_\substack{x\st{\ph}{\to}\ox}{\lm,\ve\dn
0}\lm\Hat\partial_\ve\ph(x),
\end{eqnarray}
where we can equivalently put $\ve=0$ if $\ph$ is l.s.c.\ around
$\ox$ and $X$ is Asplund. Note that
$\partial^\infty\ph(\ox)=\{0\}$ if $\ph$ is locally Lipschitzian
around $\ox$. The converse implication also holds provided that
$\ph$ is l.s.c.\ around $\ox$, that $X$ is Asplund, and that $\ph$
satisfies the so-called ``sequential normal epi-compactness''
property at $\ox$ (see below), which is always the case when $X$
is finite-dimensional.

Given a set $\O\subset X$ with its indicator function
$\dd(\cdot;\O)$ defined by $\delta(x;\O):=0$ for $x\in\O$ and by
$\delta(x;\O):=\infty$ otherwise, we construct the {\em
Fr\'echet/regular normal cone} and {\em
limiting/basic/Mordukhovich normal cone} to $\O$ at $\ox\in\O$ by,
respectively,
\begin{equation}\label{2.4}
\Hat N(\ox;\O):=\Hat\partial\delta(\ox;\O)\quad\textrm{and}\quad
N(\ox;\O):=\partial\delta(\ox;\O)
\end{equation}
via the corresponding subdifferential of the indicator function.
If follows from (\ref{2.4}) that $\Hat N(\ox;\O)\subset N(\ox;
\O)$. A set $\O$ is {\em normally regular} at $\ox$ if $\Hat
N(\ox;\O)= N(\ox;\O)$; the latter is the case of convex and some
other ``nice" sets.

Recall further that $\O$ is {\em sequentially normally compact}
(SNC) at $\ox\in\O$ if for any sequences $\ve_n\dn 0$,
$x_n\st{\O}\to\ox$, and $x^*_n\in\Hat N_{\ve_n}(x_n;\O):=\Hat \partial_{\ve_n} \delta(\ox;\O)$ we have
\begin{eqnarray*}
\big[x^*_n\st{w^*}\to 0\big]\Longrightarrow\big[\|x^*_n\|\to
0\big]\;\mbox{ as }\;n\to\infty,
\end{eqnarray*}
where $\ve_n$ can be omitted if $\O$ is locally closed around
$\ox$ and the space $X$ is  Asplund. A function $\varphi:X\to
\Bar\R$ is {\em sequentially normally epi-compact} (SNEC) at a
point $\ox\in \dom \varphi$ if its epigraph is SNC at
$(\ox,\varphi(\ox))$. Besides the finite dimensionality, the
latter properties hold under certain Lipschitzian behavior; see,
e.g., \cite[Subsections~1.1.4 and 1.2.5]{M1}.

Having an arbitrary (possibly infinite and noncompact) index set
$T$ as in \eqref{1.1}, we consider the product space of
multipliers $\R^T:=\{\lm=(\lm_t)|\;t\in T\}$ with $\lm_t\in\R$ for
$t\in T$ and denote by $\Tilde\R^T$ the collection of $\lm\in\R^T$
such that $\lm_t\ne 0$ for finitely many $t\in T$. The {\em
positive cone} in $\Tilde\R^T$ is defined by
\begin{eqnarray}\label{2.6}
\tR:=\big\{\lm\in\Tilde\R^T\big|\;\lm_t\ge 0\;\mbox{ for all }\;t\in
T\big\}.
\end{eqnarray}

\section{Qualification Conditions for Infinite Constraint Systems}
\setcounter{equation}{0}

This section is devoted to studying the set of {\em feasible
solutions} to the original optimization problem \eqref{1.1}
defined by the infinite constraint systems of inequalities and
equalities
\begin{eqnarray}\label{3.1}
\left\{\begin{array}{ll}
g_t(x)\le 0,\;\;t\in T,\\
h(x)=0,
\end{array}\right.
\end{eqnarray}
where $T$ is an arbitrary index set, and where the functions
$g_t:X\to\Bar\R$, $t\in T$, and the mapping $h:X\to Y$ are
differentiable but may not be linear and/or convex. As in
(\ref{1.2}), the set of feasible solutions to \eqref{1.1}, i.e.,
those $x\in X$ satisfying \eqref{3.1}, is denoted by
$\O$.\vspace*{0.05in}

Our {\em standing assumptions} throughout the paper (unless otherwise stated) are as
follows:\vspace*{0.05in}\\\vspace*{0.1in} {\bf (SA)} {\em For any
$\ox\in\O$ the functions $g_t$, $t\in T$, are Fr\'echet
differentiable at $\ox$ and the mapping $h$ is strictly
differentiable at $\ox$. The set $\{\nabla g_t(\ox)|\;t\in T\}$ is
bounded in $X^*$}.\vspace*{0.1in}

Recall that a mapping $h\colon X\to Y$ is {\em strictly
differentiable} at $\ox$ with the (strict) derivative $\nabla
h(\ox)\colon X\to Y$ if
\begin{equation*}
\lim_{x,x'\to\ox}\frac{h(x)-h(x')-\nabla
h(\ox)(x-x')}{\|x-x'\|}=0.
\end{equation*}
The latter holds automatically when $h$ is continuously
differentiable around $\ox$.

In addition to the standing assumptions (SA), we often impose some
stronger requirements on the inequality constraint functions $g_t$
that postulate a certain uniformity of their behavior with respect
to the index parameter $t\in T$. We say that the functions
$\{g_t\}_{t\in T}$ are {\em uniformly Fr\'echet differentiable} at
$\ox$ if
\begin{eqnarray}\label{3.3}
s(\eta):=\sup_{t\in T}\sup_\substack{x\in\B_{\eta}(\ox)}{x\ne
\ox}\frac{|g_t(x)-g_t(\ox)-\la\nabla
g_t(\ox),x-\ox\ra|}{\|x-\ox\|}\to 0\;\mbox{ as }\;\eta\dn 0.
\end{eqnarray}
Similarly, the functions $\{g_t\}_{t\in T}$ are {\em uniformly
strictly differentiable} at $\ox$ if condition (\ref{3.3}) above
is replaced by a stronger one:
\begin{eqnarray} \label{3.4}
r(\eta):=\sup_{t\in T}\sup_\substack{x,x'\in\B_{\eta}(\ox)}{x\ne
x'}\frac{|g_t(x)-g_t(x')-\la\nabla
g_t(\ox),x-x'\ra|}{\|x-x'\|}\to 0\;\mbox{ as }\;\eta\dn 0,
\end{eqnarray}
which clearly implies the strict differentiability of each
function $g_t$, $t\in T$, at $\ox$.\vspace*{0.05in}

Let us present some sufficient conditions ensuring the fulfillment
of all the assumptions formulated above for infinite families of
inequality constraint functions.

\begin{Proposition}\label{p3} {\bf (compact index sets).} Let $T$ be a
compact metric space, let the functions $g_t$ in \eqref{3.1} be
Fr\'echet differentiable around $\ox$ for each $t\in T$, and let
the mapping $(x,t)\in X\times T\mapsto\nabla g_t(x)\in X^*$ be
continuous on $\B_\eta(\ox)\times T$ for some $\eta>0$. Then the
standing assumptions {\rm (SA)} as well as {\rm(\ref{3.3})} and
{\rm(\ref{3.4})} are satisfied.
\end{Proposition}
{\bf Proof.} It is easy to see that our standing assumptions (SA)
hold, since $\|\nabla g_t(\ox)\|$ is assumed to be continuous on
the compact space $T$  being hence bounded. It suffices to prove
that (\ref{3.4}) holds, which surely implies \eqref{3.3}.

Arguing by contradiction, suppose that (\ref{3.4}) fails. Then
there are $\ve>0$, sequences $\{t_n\}\subset T$, $\{\eta_n\}\dn
0$, and $\{x_n\}$, $\{x^\prime_n\}\subset\B_{\eta_n}(\ox)$ such
that
\begin{equation}\label{3.5}
\frac{|g_{t_n}(x_n)-g_{t_n}(x'_n)-\la \nabla
g_{t_n}(\ox),x_n-x'_n\ra|}{\|x_n-x'_n\|}\ge
\ve-\frac{1}{n}\;\mbox{ for all large}\;n\in\IN.
\end{equation}
Since $T$ is a compact metric space, there is a subsequence of
$\{t_n\}$ converging (without relabeling) to some $\bar t\in T$.
Applying the classical Mean Value Theorem to (\ref{3.5}), we find
$\theta_n\in[x_n,x^\prime_n]:=\co\{x_n,x^\prime_n\}$ such that
\begin{eqnarray*}
\frac{\ve}{2}&<&\frac{|\la \nabla
g_{t_n}(\theta_n),x_n-x'_n\ra-\la
\nabla g_{t_n}(\ox),x_n-x'_n\ra|}{\|x_n-x'_n\|}\le\|\nabla g_{t_n}(\theta_n)-\nabla g_{t_n}(\ox)\|\\
&<&\|\nabla g_{t_n}(\theta_n)-\nabla g_{\bar t}(\ox)\|+\|\nabla g_{\bar t}(\ox)-\nabla g_{t_n}(\ox)\|
\end{eqnarray*}
for all large $n\in\IN$. This contradicts the continuity of the
mapping $(x,t)\in X\times T\mapsto \nabla g_t(x)$ on
$\B_\eta(\ox)\times T$ and thus completes the proof of the
proposition. $\h$\vspace*{0.1in}

Next we recall a well-recognized constraint qualification
condition, which is often used in problems of nonlinear and
nonconvex semi-infinite programming.

\begin{Definition}\label{cd1} {\bf (Extended Mangasarian-Fromovitz Constraint Qualification).}
The infinite system {\rm(\ref{3.1})} satisfies the {\sc Extended
Mangasarian-Fromovitz Constraint Qualification (EMFCQ)} at $\ox\in
\O$ if the derivative operator $\nabla h(\ox)\colon X\to Y$ is
surjective and if there is $\Tilde x\in X$ such that $\nabla
h(\ox)\tilde x=0$ and that
\begin{eqnarray}\label{3.6}
\la\nabla g_t(\ox),\tilde x\ra<0\;\mbox{ for all }\;t\in
T(\ox):=\big\{t\in T\big|\;g_t(\ox)=0\big\}.
\end{eqnarray}
\end{Definition}

It is clear that in the case of a finite index set $T$ and a
finite-dimensional space $Y$ in \eqref{3.1} the EMFCQ condition
reduced to the classical Mangasarian-Fromovitz Constraint
Qualification (MFCQ) in nonlinear programming. In the case of SIP
problems the EMFCQ was first introduced in \cite{JTW} and then
extensively studied and applied in semi-infinite frameworks with
$X=\R^m$ and $Y=\R^n$; see, e.g., \cite{HK,KH,LS,Sh}, where the
reader can find its relationships with other constraint
qualifications for SIP problems.

To the best of our knowledge, the vast majority of nonconvex
semi-infinite programs are usually considered with the general
assumptions that the index set $T$ is compact, the functions $g_t$
are continuously differentiable, and the mapping $(x,t)\mapsto
\nabla g_t(x)$ is continuous on $X\times T$. Under these
assumptions and the EMFCQ formulated above, several authors derive
the Karush-Kuhn-Tucker (KKT) necessary optimality conditions of
the following type: If $\ox$ is an optimal solution to (\ref{1.1})
with $f\in{\cal C}^1$ and $h=(h_1, h_2,\ldots,h_n)$, then there
are $\lm\in\tR$ from (\ref{2.6}) and $\mu\in \R^n$ such that
\begin{eqnarray}\label{kkt}
0=\nabla f(\ox)+\sum_{t\in T(\ox)}\lm_t\nabla
g_t(\ox)+\sum_{j=1}^n \mu_j\nabla h_j(\ox).
\end{eqnarray}

We are not familiar with any results in the literature on
nonconvex infinite programming that apply to problems with
noncompact index sets $T$. The following example shows that the
KKT optimality conditions in form \eqref{kkt} may fail for
nonconvex SIP with countable constraints even under the
fulfillment of the EMFCQ.

\begin{Example}\label{ex1} {\bf (violation of KKT for nonconvex SIP with countable
\index sets under EMFCQ).}  {\rm Consider problem (\ref{1.1}) with
countable inequality constraints given by
\begin{eqnarray}\label{3.7}
\left\{\begin{array}{ll}
\mbox{minimize }\;(x_1+1)^2+x_2\;\mbox{ subject to}\\
x_1+1\le 0\;\mbox{ and }\;\disp\frac{1}{3n}x_1^3-x_2\le
0\;\mbox{ for all }\;n\in\IN\setminus\{1\}\;\mbox{ with }\;
(x_1,x_2)\in\R^2.
\end{array}\right.
\end{eqnarray}
Let $X:=\R^2$, $Y:=\{0\}$, $f(x_1,x_2):=(x_1+1)^2+x_2$, $T:=\IN$,
$g_1(x_1,x_2):=x_1+1$, and
$g_n(x_1,x_2):=\disp\frac{1}{3n}x^3_1-x_2$ for all $n\in
\IN\setminus\{1\}$. Observe that $\ox:=(-1,0)$ is a global
minimizer for problem (\ref{3.7}) and that $T(\ox)=\{1\}$ for the
active index set in \eqref{3.6} . It is easy to check that the
EMFCQ holds at $\ox$ while there is no Lagrange multiplier
$\lm\in\R_+$ satisfying the KKT optimality condition \eqref{kkt}
at $\ox$. Indeed, we have $\la \nabla g_1(\ox),\;(-1,0)\ra=-1<0$,
and the following equation does not admit any solution for $\lm\ge
0$:
\begin{equation*}
(0,0)=\nabla f(\ox)+\lm\nabla g_1(\ox)=(0,1)+(\lm,0).
\end{equation*}}
\end{Example}\vspace*{0.05in}

Now we introduce a new extension of the MFCQ condition to the
infinite programs under consideration, which plays a crucial role
throughout the paper.

\begin{Definition}\label{cd2} {\bf (Perturbed Mangasarian-Fromovitz
Constraint Qualification).} We say that the infinite system {\rm
(\ref{3.1})} satisfies the {\sc Perturbed Mangasarian-Fromovitz
Constraint Qualification (PMFCQ)} at $\ox\in\O$ if the derivative
operator $\nabla h(\ox)\colon X\to Y$ is surjective and if there
is $\tilde x\in X$ such that $\nabla h(\ox)\tilde x=0$ and that
\begin{eqnarray}\label{3.8}
\inf_{\ve>0}\sup_{t\in T_\ve(\ox)}\la \nabla g_t(\ox),\tilde
x\ra<0\;\mbox{ with }\;T_\ve(\ox):=\big\{t\in T\big|\;g_t(\ox)\ge
-\ve\big\}.
\end{eqnarray}
\end{Definition}

In contrast to the EMFCQ, the active index set in (\ref{3.8}) is
{\em perturbed} by a small $\ve>0$. Since $T(\ox)\subset
T_\ve(\ox)$ for all $\ve>0$, the PMFCQ is stronger than the EMFCQ.
However, as shown in Section~4 and Section~5, the new condition is
much more appropriate for applications to semi-infinite and
infinite programs with general (including compact) index sets than
the EMFCQ.\vspace*{0.05in}

The following proposition reveals some assumptions on the initial
data of \eqref{3.1} ensuring the equivalence between the PMFCQ and
EMFCQ.

\begin{Proposition}\label{p4} {\bf (PMFCQ from EMFCQ).} Let $T$ be a
compact metric space, and let $\ox\in\O$ in \eqref{3.1}. Assume
that the function $t\in T\mapsto g_t(\ox)$ is upper semicontinuous
$($u.s.c.$)$ on $T$, that the derivative mapping $\nabla
h(\ox)\colon X\to Y$ is surjective, and that there is $\tilde x\in
X$ with the following properties: $\nabla h(\ox)\tilde x=0$, the
function $t\in T\mapsto\la\nabla g_t(\ox),\tilde x\ra$ is u.s.c.,
and $\la\nabla g_t(\ox),\tilde x\ra<0$ for all $t\in T(\ox)$. Then
the PMFCQ condition holds at $\ox$, being thus equivalent to the
EMFCQ condition at this point.
\end{Proposition}
{\bf Proof.} Arguing by contradiction, suppose that the PMFCQ
fails at $\ox$. Then it follows from (\ref{3.8}) that there exist
sequences $\{\ve_n\}\dn 0$ and $\{t_n\}\subset T$ such that
$t_n\in T_{\ve_n}(\ox)$ and
\begin{eqnarray*}
\la\nabla g_{t_n}(\ox),\tilde x\ra\ge-\frac{1}{n}\;\mbox{ for all
}\;n\in\IN.
\end{eqnarray*}
Since $T$ is a compact metric space, we find a subsequence of
$\{t_n\}$ (no relabeling), which converges to some $\bar t\in T$.
Observe from the continuity assumptions made imply that
\begin{eqnarray*}
g_{\bar
t}(\ox)\ge\limsup_{n\to\infty}g_{t_n}(\ox)\ge\limsup_{n\to\infty}-\ve_n=0\;\mbox{
 and}\\ \la \nabla g_{\bar t}(\ox),\tilde x\ra\ge\limsup_{n\to\infty}
\la\nabla g_{t_n}(\ox),\tilde x\ra\ge
\limsup_{n\to\infty}-\frac{1}{n}=0.
\end{eqnarray*}
Thus we have that $\bar t\in T(\ox)$ and $\la\nabla g_{\bar
t}(\ox),\tilde x\ra\ge 0$, which is a contradiction that completes
the proof of the proposition. $\h$\vspace*{0.05in}

The following example shows that the EMFCQ does not imply the
PMFCQ (while not ensuring in this case the validity of the
required necessary optimality conditions as will be seen in
Sections~4 and 5) even for simple frameworks of  nonconvex
semi-infinite programs with {\em compact} index sets.

\begin{Example}\label{exa1} {\bf (EMFCQ does not imply PMFCQ for
semi-infinite programs with compact index sets).} {\rm Let
$X=\R^2$ and $T=[0,1]$ in \eqref{3.1} with $h=0$ and
\begin{eqnarray*}
g_0(x):=x_1+1\le 0,\quad g_t(x):=tx_1-x^3_2\le 0\;\mbox{ for
}\;t\in T\setminus\{0\}.
\end{eqnarray*}
It is easy to check that the functions $g_t$, $t\in T$, satisfy
our standing assumptions and that they are strictly uniformly
differentiable at the feasible point $\ox=(-1,0)$. Observe
furthermore that $T(\ox)=\{0\}$, that $T_\ve(\ox)=[0,\ve]$ for all
$\ve\in(0,1)$, and that the EMFCQ holds at $\ox$. However, for any
$d=(d_1,d_2)\in \R^2$ we have
\begin{eqnarray*}
\disp\inf_{\ve>0}\sup_{t\in T_{\ve(\ox)}}\la\nabla
g_t(\ox),d\ra&=&\inf_{\ve>0}
\sup\Big\{\la \nabla g_0(\ox),d\ra,\sup\big\{\la\nabla g_t(\ox),d\ra\big|\;t\in(0,\ve]\big\}\Big\}\\
&=&\inf_{\ve>0}\sup\Big\{d_1,\sup\{td_1\big|\;t\in
(0,\ve]\}\Big\}\ge 0,
\end{eqnarray*}
which shows that the PMFCQ does not satisfy at $\ox$. Note that
the u.s.c.\ assumption with respect of $t$ in
Propositions~\ref{p4} does not hold in this example.}
\end{Example}

It is well known in the classical nonlinear programming (when the
index set $T$ in \eqref{3.1} is finite), that the MFCQ condition
is equivalent to the Slater condition provided that all the
functions $g_t$ are convex and differentiable and that $h$ is a
linear operator. The next proposition shows that a similar
equivalence holds in the semi-infinite and infinite programming
frameworks with replacing the MFCQ by our new PMFCQ condition and
replacing the Slater by its strong counterpart well recognized in
the SIP community; see, e.g., \cite{GL} and \cite{CLMP1} for more
references and discussions.

\begin{Proposition}\label{ss} {\bf (equivalence between PMFCQ and SSC
for differentiable convex systems).} Assume that in \eqref{3.1}
all the functions $g_t$, $t\in T$, are convex and uniformly
Fr\'echet differentiable at $\ox$ and that $h=A$ is a surjective
continuous linear operator. Then the PMFCQ condition is equivalent
to the following strong Slater condition $($SSC$)$: there is $\hat
x\in X$ such that $A\hat x=0$ and
\begin{eqnarray}\label{3.9}
\sup_{t\in T} g_t(\hat x)<0.
\end{eqnarray}
\end{Proposition}
{\bf Proof.} Suppose first that the SSC holds at $\ox$, i.e.,
there are $\hat x\in X$ and $\delta>0$ such that $A\hat x=0$ and
$g_t(\hat x)<-2\delta$ for all $t\in T$. By the assumptions made
this implies that for each $\ve\in (0,\delta)$ and $t\in
T_\ve(\ox)$ we have
$$
\la\nabla g_t(\ox),\hat x-\ox\ra\le g_t(\hat x)-g_t(\ox)\le-2\delta+\ve\le-\delta.
$$
Define further $\tilde x:=\hat x-\ox$ and get $A\tilde x=A\hat
x-A\ox=0$ with $\la\nabla g_t(\ox),\tilde x\ra\le-\delta$ for all
$t\in T_\ve(\ox)$ and  $\ve\in(0,\delta)$. This clearly implies
the PMFCQ condition at $\ox$.

Conversely, assume that the PMFCQ condition holds at $\ox$. Then
there are $\ve,\eta>0$ and $\Tilde x\in X$ such that $\la\nabla
g_t(\ox),\Tilde x\ra\le-\eta$ for all $t\in T_\ve(\ox)$ and that
$A\Tilde x=0$. It follows from the assumed uniform Fr\'echet
differentiability \eqref{3.3} of $g_t$ at $\ox$ that for each
$\lm>0$ we have
\begin{equation}\label{t}
g_t(\ox+\lm\Tilde x)\le g_t(\ox)+\lm\la\nabla g_t(\ox),\Tilde
x\ra+\lm\|\Tilde x\|s\big(\lm\|\Tilde x\|\big),
\end{equation}
which readily implies that $g_t(\ox+\lm\Tilde
x)\le\lm\big(-\eta+\|\Tilde x\|s(\lm\|\Tilde x\|)\big)$ for all $t\in T_\ve(\ox)$. For
$t\notin T_\ve(\ox)$ we observe from \eqref{t} that
\begin{equation*}
g_t(\ox+\lm\Tilde x)\le-\ve+\lm\sup_{\tau\in T}\|\nabla
g_\tau(\ox)\|\cdot\|\Tilde x\|+\lm \|\Tilde x\|s\big(\lm\|\Tilde x\|\big),
\end{equation*}
which gives, combining with the above, that
\begin{equation*}
\sup_{t\in T} g_t(\ox+\lm\Tilde x)\le\max\big\{\lm\big(-\eta+\|\Tilde x\|s(\lm\|\Tilde x\|)\big),-\ve+\lm\|\Tilde x\|\big(\sup_{\tau\in T}\|\nabla g_\tau(\ox)\|+s(\lm\|\Tilde x\|)\big)\big\}.
\end{equation*}
The latter implies the existence of $\lm_0>0$ sufficiently small
such that $\sup_{t\in T}g_t(\Hat x)<0$ with $\Hat
x:=\ox+\lm_0\Tilde x$. Furthermore, it is easy to see that $A\Hat
x=A\ox+\lm_0 A\Tilde x=0$. This concludes that the SSC holds at
$\Hat x$ and thus completes the proof of the proposition.
$\h$\vspace*{0.1in}

Next we introduce another qualification condition of the {\em
closedness/Farkas-Minkowski type} for infinite inequality
constraints in \eqref{1.1}.

\begin{Definition}\label{cd3} {\bf (Nonlinear Farkas-Minkowski Constraint Qualification).}
We say that system {\rm(\ref{3.1})} with $h(x)=0$ satisfies the
{\sc Nonlinear Farkas-Minkowski constraint qualification (NFMCQ)}
at $\ox$ if the set
\begin{eqnarray}\label{fm}
\cone\big\{\big(\nabla g_t(\ox),\la\nabla
g_t(\ox),\ox\ra-g_t(\ox)\big)\big|\;t\in T\big\}
\end{eqnarray}
is weak$^*$ closed in the product space $X^*\times\R$.
\end{Definition}
In the linear case of $g_t(x)=\la a^*_t,x\ra-b_t$ for some
$(a^*_t,b_t)\in X^*\times \R$, $t\in T$, the NFMCQ condition above
reduces to the classical Farkas-Minkowski qualification condition
meaning that the set $\cone\{(a^*_t,b_t)|\; t\in T\}$ is weak$^*$
closed in $X^*\times\R$. It is well recognized that the latter
condition plays an important role in linear semi-infinite and
infinite optimization; see, e.g.,
\cite{BGW,CLMP2,DGL,DMN1,DMN2,GL} for more details and references.
Observe that the NFMCQ condition can be represented in the
following equivalent form: the set
\begin{eqnarray*}
\cone\big\{\big(\nabla g_t(\ox),g_t(\ox)\big)\big|\;t\in
T\big\}\;\mbox{ is weak$^*$ closed in $X^*\times \R$}.
\end{eqnarray*}

Let us compare the new NFMCQ condition with the other
qualification conditions discussed in this section in the case of
infinite inequality constraints.

\begin{Proposition}\label{sfm} {\bf (sufficient conditions for NFMCQ).}
Consider the constraint inequality system \eqref{3.1} with $h=0$
therein. Then the NFMCQ condition is satisfied at $\ox\in\O$ in
each of the following settings:

{\bf (i)} The index $T$ is finite and the MFCQ condition holds at
$\ox$.

{\bf (ii)} $\dim X<\infty$, the set $\{(\nabla g_t(\ox),\la\nabla
g_t(\ox),\ox\ra-g_t(\ox))|\;t\in T\}$ is compact, and the PMFCQ
condition holds at $\ox$.

{\bf(iii)} The index $T$ is a compact metric space, $\dim
X<\infty$, the mappings $t\in T\mapsto g_t(\ox)$ and $t\in
T\mapsto\nabla g_t(\ox)$ are continuous, and the EMFCQ condition
holds at $\ox$.
\end{Proposition}
{\bf Proof.} Define $\tilde g_t(x):=\la \nabla
g_t(\ox),x-\ox\ra+g_t(\ox)$ for all $x\in X$. To justify (i),
suppose that $T$ is finite and that the MFCQ condition holds at
$\ox$ for the inequality system in (\ref{3.1}). It is clear that
$\tilde g_t$ also satisfy the MFCQ at $\ox$. Since the functions
$\tilde g_t$ are linear, we observe from Proposition~\ref{ss} that
there is $\hat x\in X$ such that $\tilde g_t(\hat x)=\la \nabla
g_t(\ox),\hat x-\ox\ra+g_t(\ox)< 0$ for all $t\in T$. Thus it
follows from \cite[Proposition~6.1]{DMN1} that the NFMCQ condition
holds.

Next we consider case (ii) with $X=\R^d$ therein. Suppose that the
PMFCQ condition holds at $\ox$ and that the set $\{(\nabla
g_t(\ox),\la \nabla g_t(\ox),\ox\ra-g_t(\ox))|\;t\in T\}$ is
compact in $\R^d\times\R$. Noting that the functions $\tilde g_t$
also satisfy the PMFCQ at $\ox$, we apply Proposition~\ref{ss} to
these functions and find $\hat x\in X$ such that $\nabla
h(\ox)\hat x=0$ and that
\begin{eqnarray}\label{3.11}
\sup_{t\in T} \tilde g_t(\hat x)=\sup_{t\in T} \la \nabla g_t(\ox),
\hat x-\ox\ra+g_t(\ox)<0.
\end{eqnarray}
Let us check that $(0,0)\not \in \co\{(\nabla g_t(\ox),\la\nabla
g_t(\ox),\ox\ra-g_t(\ox))|\; t\in T\}$. Indeed, otherwise ensures
the existence of $\lm\in\tR$ with $\sum_{t\in T}\lm_t=1$ such that
\begin{eqnarray*}
(0,0)=\sum_{t\in T}\lm_t\big(\nabla g_t(\ox),\la\nabla
g_t(\ox),\ox\ra-g_t(\ox)\big).
\end{eqnarray*}
Combining the latter with (\ref{3.11}) gives us that
\begin{equation*}
0=\sum_{t\in T}\lm_t\la\nabla g_t(\ox),\hat x\ra-\sum_{t\in
T}\lm_t\big(\la\nabla g_t(\ox),\ox\ra-g_t(\ox)\big)=\sum_{t\in
T}\lm_t\tilde g_t(\hat x)\le\sup_{t\in T}\tilde g_t(\hat x)<0,
\end{equation*}
which is a contradiction. Hence employing \cite[Theorem~1.4.7]{HL}
in this setting, we have that the conic hull $\cone\{(\nabla
g_t(\ox),\la \nabla g_t(\ox),\ox\ra-g_t(\ox))|\;t\in T\}$ is
closed in $\R^{d+1}$. This fully justifies (ii). Observing finally
that (iii) follows from (ii) and Proposition~\ref{p4}, we complete
the proof of the proposition. $\h$\vspace*{0.05in}

To conclude this section, let us show that the NFMCQ and PMFCQ
conditions are independent for infinite inequality systems in
finite dimensions.

\begin{Example}\label{exa2} {\bf (independence of NFMCQ and
PMFCQ).} {\rm It is easy to check that for the constraint
inequality system from Example~\ref{exa1} the NFMCQ is satisfied
at $\ox=(-1,0)$, since the corresponding conic hull
\begin{eqnarray*}
\cone\big\{(\nabla g_t(\ox),\la\nabla
g_t(\ox),\ox\ra-g_t(\ox))\big|\;t\in
T\big\}&=&\cone\Big((1,0,-1)\cup
\{(t,0,0)\big|t\in(0,1]\}\Big)\\
&=&\big\{x\in\R^3\big|\;x_1+x_3\ge 0,\;x_1\ge 0\ge
x_3,\;x_2=0\big\}
\end{eqnarray*}
is closed in $\R^3$. On the other hand, Example~\ref{exa1}
demonstrates that the PMFCQ does not hold for this system at
$\ox$.

To show that the NFMCQ does not generally follow from the PMFCQ
(and even from the EMFCQ), consider the countable system of
inequality constraints \eqref{3.7} in $\R^2$ discussed in
Example~\ref{ex1}. When $\ox=(-1,0)$, we get $T_\ve(\ox)=\{n\in
\IN\setminus\{1\}|\;n\le \frac{1}{\ve}\}\cup\{1\}$ for the the
perturbed active index set in \eqref{3.8}. It shows that the PMFCQ
(and hence the EMFCQ) hold at $\ox$. On the other hand, the conic
hull
\begin{eqnarray*}
\mbox{cone}\big\{(\nabla g_t(\ox),\la\nabla
g_t(\ox),\ox\ra-g_t(\ox))\big|\;t\in
T\big\}=\mbox{cone}\Big[(1,0,-1)\cup\Big\{\Big(\frac{1}{n},-1,
\frac{-2}{3n}\Big)\Big|\;n\in\IN\setminus\{1\}\Big\}\Big]
\end{eqnarray*}
is not closed in $\R^3$, i.e., the NFMCQ condition is not
satisfies at $\ox$.}
\end{Example}

\section{Normal Cones to Feasible Sets of Infinite Constraints}
\setcounter{equation}{0}

This section is devoted to computing both normal cones \eqref{2.4}
to the feasible solution sets \eqref{1.2} for the class of
nonconvex semi-infinite/infinite programs \eqref{1.1} under
consideration in the paper. These calculus results are certainly
of independent interest while they play a crucial role in deriving
necessary optimality conditions for \eqref{1.1} in Section~5.\vspace*{0.05in}

The first main theorem gives precise calculations of both
Fr\'echet and limiting normal cones to the set $\O$ of feasible
solutions in \eqref{1.2} under the new Perturbed
Mangasarian-Fromovitz Constraint Qualification of
Definition~\ref{cd2}. Preliminary we present a known result
from functional analysis whose simple proof is given for the reader's convenience.

\begin{Lemma}\label{cr} {\bf(weak$^*$ closed images of adjoint operators).} Let $A:X\to Y$ be a surjective
continuous linear operator. Then the image of its adjoint operator $A^*(Y^*)$ is a weak$^*$ closed
subspace of $X^*$.
\end{Lemma}
{\bf Proof.} Define $C:=A^*(Y^*)\subset X^*$ and pick any $n\in\IN$. We claim that the set $A_n:=C\cap n\B_{X^*}$
is weak$^*$ closed in $X^*$. Considering a net $\{x^*_\nu\}_{\nu\in\N}\subset A_n$ weak$^*$ converging to
$x^*\in X^*$ and taking into account that the ball $\B_{X^*}$ is weak$^*$ compact in $X^*$, we get $x^*\in n\B_{X^*}$.
By construction there is a net $\{y^*_\nu\}_{\nu\in \N}\subset Y^*$ satisfying $x^*_\nu=A^*y^*_\nu$ whenever $\nu\in\N$.
It follows from the surjectivity of $A$ that
$$
\|x^*_\nu\|=\|A^*y^*_\nu\|\ge\kappa\|y^*_\nu\|\;\mbox{ for all }\;\nu\in\N,
$$
where $\kappa:=\inf\{\|A^*y^*\|\;\mbox{over}\;\|y^*\|=1\}\in(0,\infty)$; see, e.g., \cite[Lemma~1.18]{M1}. Hence
$\|y^*_\nu\|\le n\kappa^{-1}$ for all $\nu\in\N$. By passing to a subnet, suppose that $y^*_\nu$ weak$^*$ converges to some
$y^*\in Y^*$ for which $x^*=A^*y^*\in A_n$. Thus we have that the set $A_n=C\cap n\B_{X^*}$ is weak$^*$ closed for all
$n\in\IN$. The classical Banach-Dieudonn\'e-Krein-\v Smulian theorem yields therefore that the set $C$ is weak$^*$ closed in
$X^*$. $\h$\vspace*{0.05in}

Now we are ready to establish the main result of this section.

\begin{Theorem}\label{thm41} {\bf (Fr\'echet and limiting normals to infinite
constraint systems).} Let $\ox\in\O$ for the set of feasible
solutions \eqref{1.2} to the infinite system \eqref{3.1}
satisfying the PMFCQ at $\ox$. Assume in addition that the
inequality constraint functions $g_t$, $t\in T$, are uniformly
Fr\'echet differentiable at $\ox$. Then the Fr\'echet normal cone
to $\O$ at $\ox$ is computed by
\begin{eqnarray}\label{4.1}
\Hat N(\ox;\O)=\disp\bigcap_{\ve>0}\cl\cone\big\{\nabla
g_t(\ox)\big|\;t\in T_\ve(\ox)\big\}+\nabla h(\ox)^*(Y^*).
\end{eqnarray}
If furthermore the functions $g_t$, $t\in T$, are uniformly
strictly differentiable at $\ox$, then the limiting normal cone to
$\O$ at $\ox$ is also computed by
\begin{eqnarray}\label{4.2}
N(\ox;\O)=\disp\bigcap_{\ve>0}\cl\cone\big\{\nabla
g_t(\ox)\big|\;t\in T_\ve(\ox)\big\}+\nabla h(\ox)^*(Y^*),
\end{eqnarray}
and thus the set $\O$ of feasible solutions is normally regular at
$\ox$.
\end{Theorem}
{\bf Proof.} First we justify \eqref{4.1} under the assumptions
made. It follows from the PMFCQ and the uniform Fr\'echet
differentiability of $g_t$ at $\ox$ that there are $\tilde\ve>0$,
$\delta>0$, and $\tilde x\in X$ such that $\nabla h(\ox)\tilde
x=0$ and
\begin{eqnarray}\label{4.3}
\disp\sup_{t\in T_\ve(\ox)}\la \nabla g_t(\ox),\tilde
x\ra<-\delta\quad {\rm for\ all}\ \ve\leq \tilde\ve.
\end{eqnarray}
Let us prove the inclusion ``$\supset$'' in (\ref{4.1}). To
proceed, fix any $\ve\in(0,\tilde \ve)$ and pick an arbitrary
element $x^*$ belonging to the right-hand side of (\ref{4.1}).
Then there exist a net  $(\lm_\nu)_{\nu\in \N}\subset\tR$ and a
dual element $y^*\in Y^*$ satisfying
\begin{eqnarray}\label{4.4}
x^*=w^*-\lim_\nu\sum_{t\in T_\ve(\ox)} \lm_{t\nu}\nabla g_t(\ox)+\nabla h(\ox)^*y^*.
\end{eqnarray}
Combining the latter with (\ref{4.3}) gives us
\begin{eqnarray}\begin{array}{ll}\label{4.5}
\la x^*,\tilde x\ra&=\disp\lim_\nu\sum_{t\in
T_\ve(\ox)}\lm_{t\nu}\la \nabla g_t(\ox),
\tilde x\ra+\la\nabla h(\ox)^*y^*,\tilde x\ra\\
&\disp\le\liminf_\nu\sum_{t\in T_\ve(\ox)} \lm_{t\nu}(-\delta)+\la
y^*,\nabla h(\ox)\tilde x\ra= -\delta\limsup_\nu\sum_{t\in
T_\ve(\ox)}\lm_{t\nu}.
\end{array}
\end{eqnarray}
It follows further that for each $\eta>0$ and $x\in\O\cap
\B_\eta(\ox)$ we have
\begin{eqnarray*}
\la x^*,x-\ox\ra&=& \lim_\nu\sum_{t\in T_\ve(\ox)}\lm_{t\nu}\la
\nabla g_t(\ox),x-\ox\ra+\la\nabla h(\ox)^*y^*,x-\ox\ra\\
&\le&\limsup_\nu\sum_{t\in T_\ve(\ox)}\lm_{t\nu}\Big(
g_t(x)-g_t(\ox)+\|x-\ox\|s(\eta)\Big)+
\la y^*,\nabla h(\ox)(x-\ox)\ra\\
&\le&\limsup_\nu\sum_{t\in T_\ve(\ox)}
\lm_{t\nu}\Big(\ve+\|x-\ox\|s(\eta)\Big)+ \|y^*\|\Big
(\|h(x)-h(\ox)\|+o(\|x-\ox\|)\Big)\\
&\le&\Big(\ve+\|x-\ox\|s(\eta)\Big)\limsup_\nu\sum_{t\in
T_\ve(\ox)}\lm_{t\nu}+\|y^*\|o(\|x-\ox\|).
\end{eqnarray*}
Taking now the estimate (\ref{4.5}) into account implies that
\begin{eqnarray*}
\la x^*,x-\ox\ra\le-\frac{\la x^*,\tilde
x\ra}{\delta}\Big(\ve+\|x-\ox\|s(\eta)\Big)+o(\|x-\ox\|)\|y^*\|,
\end{eqnarray*}
which yields in turn by $\ve\dn 0$ that
\begin{eqnarray*}
\la x^*,x-\ox\ra\le-\frac{\la x^*,\tilde
x\ra}{\delta}\|x-\ox\|s(\eta)+o(\|x-\ox\|)\|y^*\|.
\end{eqnarray*}
Since $s(\eta)\dn 0$ as $\eta\dn 0$, it follows from the latter inequality that
\begin{eqnarray*}
\disp\limsup_{x\st{\O}\to\ox}\frac{\la x^*,
x-\ox\ra}{\|x-\ox\|}\le 0,
\end{eqnarray*}
which means that $x^*\in\Hat N(\ox;\O)$ and thus justifies the
inclusion ``$\supset$'' in (\ref{4.1}).

Next we prove the inclusion ``$\subset$'' in (\ref{4.2}) under the
assumption that $g_t$ are uniformly strictly differentiable at
$\ox$. This immediately implies the inclusion ``$\subset$'' in
(\ref{4.1}) under the latter assumption, while we note that
similar arguments justify the inclusion ``$\subset$'' in
(\ref{4.1}) under merely the uniform Fr\'echet differentiability
of $g_t$ at $\ox$.

To proceed with proving the inclusion ``$\subset$'' in
(\ref{4.2}), define the set
\begin{equation}\label{A}
A_\ve:=\cl{\rm cone}\big\{\nabla g_t(\ox)\big|\;t\in
T_\ve(\ox)\big\}+\nabla h(\ox)^*(Y^*)\;\mbox{ for }\;\ve>0.
\end{equation}
 Arguing by contradiction, pick an arbitrary element $x^*\in N(\ox;\O)\setminus\{0\}$ and suppose that $
 x^*\notin A_\ve$ for some $\ve\in(0,\tilde\ve)$. We first claim that the set $A_\ve$ is weak$^*$ closed in $X^*$ for all $\ve\leq \tilde\ve$. To justify, take  an arbitrary net $(u_\nu^*)_{\nu\in \N}\subset A_\ve$ weak$^*$ converging to some $u^*\in X^*$. Hence there are nets $(\lm_\nu)_{\nu\in \N}\subset \tR$, $(y^*_\nu)_{\nu\in \N} \subset Y^*$ such that
\begin{equation*}
u^*_\nu=\sum_{t\in T_\ve(\ox)}\lm_{t\nu}\nabla g_t(\ox)+\nabla h(\ox)^*y^*_\nu\quad \st{w^*}\to u^*.
\end{equation*}
Similarly to the proof of (\ref{4.5}) we derive the inequality
\begin{eqnarray*}
\la u^*,\tilde x\ra\leq-\delta\limsup_\nu\sum_{t\in T_\ve(\ox)} \lm_{t\nu}.
\end{eqnarray*}
Moreover, we have
\[
\|u^*_\nu-\nabla h(\ox)^*y^*_\nu\|=\|\sum_{t\in T_\ve(\ox)}\lm_{t\nu}\nabla g_t(\ox)\|\leq \sup_{\tau\in T_\ve(\ox)}\|\nabla g_\tau(\ox)\|\sum_{t\in T_\ve(\ox)}\lm_{t\nu}.
\]
It follows from two inequalities above that the net $\{u^*_\nu-\nabla h(\ox)^*y^*_\nu\}_{\nu\in \N}$ is bounded in $X^*$. By the classical Alaoglu-Bourbaki theorem, there is a subnet of $\{u^*_\nu-\nabla h(\ox)^*y^*_\nu\}$ (without relabeling) weak$^*$ converging to some $v^*\in \cl{\rm cone}\big\{\nabla g_t(\ox)\big|\;t\in
T_\ve(\ox)\big\}$. Thus the net $\{\nabla h(\ox)^*y^*_\nu\}$ weak$^*$ converges to $u^*-v^*$. Due to Lemma \ref{cr}, there is $y^*\in Y^*$ such that $u^*-v^*=\nabla h(\ox)^*y^*$. This implies that $u^*=v^*+\nabla h(\ox)^*y^*\in A_\ve$ and ensures that $A_\ve$ is weak$^*$ closed in $X^*$. Since $x^*\notin A_\ve$, we conclude from the classical separation theorem that there are $x_0\in X$ and $c>0$ satisfying
\begin{eqnarray}\label{4.6}
\la x^*,x_0\ra\ge 2c>0\ge\la\nabla g_t(\ox),x_0\ra+\la\nabla h(\ox)^*y^*,x_0\ra=\la\nabla g_t(\ox),x_0\ra+\la y^*,\nabla h(\ox)x_0\ra
\end{eqnarray}
for all $t\in T_\ve(\ox)$ and $y^*\in Y^*$; hence $\nabla h(\ox)x_0=0$. Define further
$$
\disp\hat x:=x_0+\frac{c}{\|x^*\|\cdot\|\tilde x\|}
\tilde x
$$
and observe that $\nabla h(\ox)\hat x=0$. Moreover, it follows from (\ref{4.6}) and the PMFCQ that
\begin{eqnarray}\label{4.7}
\la x^*,\hat x\ra=\disp\la x^*,x_0+\frac{c}{\|x^*\|\cdot\|\tilde x\|}\tilde x\ra\ge 2c+\frac{c}{\|x^*\|\cdot\|\tilde x\|}\la x^*,
\tilde x\ra\ge 2c-c=c\;\mbox{ and}
\end{eqnarray}
\begin{eqnarray}\label{4.8}
\la\nabla g_t(\ox),\hat x\ra=\disp\la\nabla g_t(\ox),x_0\ra+\frac{c}{\|x^*\|\cdot\|\tilde x\|}\la\nabla g_t(\ox),
\tilde x\ra\le-\disp\frac{\delta c}{\|x^*\|\cdot\|\tilde x\|}=-\tilde \delta
\end{eqnarray}
for all $t\in T_\ve(\ox)$ with $\tilde\delta:=\disp\frac{\delta c}{\|x^*\|\cdot\|\tilde x\|}>0$. Observing that $\hat x\ne 0$
by (\ref{4.8}), suppose without loss of generality that $\|\hat x\|=1$. Furthermore, we get from definition of the limiting
normal cone that there are sequences $\ve_n\dn 0$, $\eta_n\dn 0$, $x_n\st{\O}\to \ox$, and $x^*_n\st{w^*}\to x^*$ as $n\to\infty$ with
\begin{eqnarray}\label{4.9}
\la x^*_n,x-x_n\ra\le\ve_n\|x-x_n\|\;\mbox{ for all }\;x\in\B_{\eta_n}(x_n)\cap\O,\quad n\in\IN.
\end{eqnarray}
Since the mapping $h$ is strictly differentiable at $\ox$ with the surjective derivative $\nabla h(\ox)$, it follows from the
Lyusternik-Graves theorem  (see, e.g., \cite[Theorem~1.57]{M1}) that $h$ is {\em metrically regular} around $\ox$, i.e., there
are neighborhoods $U$ of $\ox$ and $V$ of $0=h(\ox)$ and a constant $\mu>0$  such that
\begin{equation}\label{mr}
{\rm dist}\big(x; h^{-1}(y)\big)\le\mu\|y-h(x)\|\;\mbox{ for any }\;x\in U\;\mbox{ and }\;y\in V.
\end{equation}
Since $h(x_n)=0$ and $\nabla h(\ox)\hat x=0$, we have
\begin{eqnarray*}
\|h(x_n+t\hat x)\|=\|h(x_n+t\hat x)-h(x_n)-\nabla h(\ox)(t\hat x)\|=o(t)\;\mbox{ for each small }\ t>0.
\end{eqnarray*}
Thus the metric regularity (\ref{mr}) implies that for any small $t>0$ there is $ x_t\in h^{-1}(0)$ with $\|x_n+t\hat x-x_t\|=o(t)$
when $x_n\in U$. This allows us to find $\tilde\eta_n<\eta_n$ and $\tilde x_n:=x_{\tilde\eta_n}\in h^{-1}(0)$ satisfying
$\tilde\eta_n+o(\tilde \eta_n)\le\eta_n$ and $\|x_n+\tilde\eta_n\hat x-\tilde x_n\|=o(\tilde\eta_n)$. Note that
$$
\|x_n-\tilde x_n\|\le\tilde\eta_n\|\hat x\|+\|x_n+\tilde\eta_n\hat x-\tilde x_n\|=\tilde\eta_n+o(\tilde\eta_n)\le\eta_n,
$$
i.e.,  $\tilde x_n\in\B_{\eta_n}(x_n)$. Observe further that
$$
\|x_n-\tilde x_n\|\ge\tilde\eta_n\|\hat x\|-\|x_n+\tilde\eta_n\hat x-\tilde x_n\|=\tilde\eta_n-o(\tilde\eta_n).
$$
By the classical uniform boundedness principle there is a constant $M$ such that $M>\|x^*_n\|$ for all $n\in\IN$ due to $x_n^*\st{w^*}\to x^*$ as $n\to\infty$. It follows from (\ref{4.7}) that $\la x^*_n,\hat x\ra>0$ for $n\in\IN$ sufficiently large. Then we have
\begin{eqnarray*}
\disp\frac{\la x^*_n,\tilde x_n-x_n\ra}{\|\tilde x_n-x_n\|}&=&\disp\frac{\la x^*_n,\tilde x_n-\tilde\eta_n\hat x-x_n\ra}{\|\tilde x_n-x_n\|}
+\frac{\la x^*_n,\tilde\eta_n \hat x\ra}{\|\tilde x_n-x_n\|}\\
&\ge&\disp-M\frac{\|\tilde x_n-\tilde\eta_n\hat x-x_n\|}{\|\tilde x_n-x_n\|}+\tilde\eta_n\frac{\la x^*_n,\hat x\ra}{\|\tilde x_n-x_n\|}\\
&\ge&\disp-M\frac{o(\tilde\eta_n)}{\tilde\eta_n-o(\tilde\eta_n)}+\frac{\tilde\eta_n}{\tilde\eta_n+o(\tilde\eta_n)}\la x_n^*,\hat x\ra.
\end{eqnarray*}
Since $o(\tilde\eta_n)/\tilde\eta_n\to 0$ when $n\to\infty$, the latter inequalities yield that
\begin{eqnarray*}
\liminf_{n\to\infty}\disp\frac{\la x^*_n,\tilde x_n-x_n\ra}{\|\tilde x_n-x_n\|}\ge\la x^*,\hat x\ra.
\end{eqnarray*}
Combining this  with (\ref{4.7}) and (\ref{4.9}) gives us that $\tilde x_n\notin\O$ for all large $n\in\IN$.

Now define $u_n:=x_n+\tilde\eta_n\hat x-\tilde x_n$ and get $\|u_n\|=o(\tilde\eta_n)$ and $\|\tilde x_n+u_n-x_n\|=\tilde\eta_n$
by the arguments above. It follows from our standing assumptions (SA), condition (\ref{3.4}), and inequality (\ref{4.8}) that for
each $t\in T_\ve(\ox)$ we have
\begin{eqnarray*}
-\disp\tilde \delta&\geq&\disp\frac{\la\nabla g_t(\ox),\tilde\eta_n\hat x\ra}{\tilde\eta_n}=\disp\frac{\la\nabla g_t(\ox),\tilde x_n+u_n-x_n\ra}{\|\tilde x_n+u_n-x_n\|}=\disp\frac{\la \nabla g_t(\ox),\tilde x_n-x_n\ra}{\|\tilde x_n+u_n-x_n\|}
+\disp\frac{\la\nabla g_t(\ox),u_n\ra}{\|\tilde x_n+u_n-x_n\|}\\
&\ge&\disp\frac{\la\nabla g_t(\ox),\tilde x_n-x_n\ra}{\|\tilde x_n-x_n\|}\frac{\|\tilde x_n-x_n\|}{\|\tilde x_n+u_n-x_n\|}+\disp
\frac{\la\nabla g_t(\ox),u_n\ra}{\|\tilde x_n+u_n-x_n\|}\\
&\ge&\disp\Big(\frac{g_t(\tilde x_n)-g_t(x_n)}{\|\tilde x_n-x_n\|}-r(\hat\eta_n)\Big)\frac{\|\tilde x_n-x_n\|}{\|\tilde x_n+u_n-x_n\|}
-\sup_{\tau\in T_\ve(\ox)}\|\nabla g_\tau(\ox)\|\frac{o(\tilde\eta_n)}{\tilde\eta_n}\\
&\geq&  \disp \Big(\frac{g_t(\tilde x_n)}{\|\tilde x_n-x_n\|}-r(\hat \eta_n)\Big)\frac{\|\tilde x_n-x_n\|}{\|\tilde x_n+u_n-x_n\|}- \sup_{\tau\in T}\|\nabla g_\tau(\ox)\|\frac{o(\tilde\eta_n)}{\tilde\eta_n},
\end{eqnarray*}
where $\hat\eta_n:=\max\{\|x_n-\ox\|$ and $\|\tilde x_n-\ox\|\}\to 0$ as $n\to\infty$. Note that
\begin{eqnarray*}
\frac{\tilde\eta_n-o(\tilde\eta_n)}{\tilde\eta_n}\le\frac{\|\tilde x_n-x_n\|}{\|\tilde x_n+u_n-x_n\|}\le\frac{\tilde\eta_n+o(\tilde \eta_n)}{\tilde\eta_n},
\end{eqnarray*}
which implies that $\disp\frac{\|\tilde x_n-x_n\|}{\|\tilde x_n+u_n-x_n\|}\to 1$ as $n\to\infty$. Furthermore, since $r(\hat\eta_n)\to 0$ and $\disp \frac{o(\tilde \eta_n)}{\tilde\eta_n}\to 0$ as $n\to \infty$, we have ${g_t(\tilde x_n)}\le-\disp\frac{\tilde\delta}{2}\|\tilde x_n-x_n\|\le 0$ for each $t\in T_\ve(\ox)$ when $n\in\IN$ is sufficiently large. Indeed, assuming otherwise that $t\notin T_\ve(\ox)$ gives us
\begin{eqnarray*}
g_t(\tilde x_n)&\le&g_t(\ox)+\la\nabla g_t(\ox),x_n-\ox\ra+\|x_n-\ox\|r(\hat\eta_n)\\
&\le&-\ve+\sup_{\tau\in T}\|\nabla g_\tau(\ox)\|\hat\eta_n+\hat\eta_n r(\hat\eta_n)\le 0\;\mbox{ for all large }\;n\in\IN.
\end{eqnarray*}
Thus $g_t(\tilde x_n)\le 0$ for all $t\in T$ and also $h(\tilde x_n)=0$ when $n\in\IN$ is sufficiently large, i.e., $\tilde x_n\in\O$,
a contradiction. Hence we conclude that $N(\ox;\O)\subset A_\ve$ for all $\ve\in(0,\tilde\ve)$, which implies the inclusion ``$\subset$'' in (\ref{4.2}) and completes the proof of the theorem. $\h$\vspace*{0.05in}

Let us show now that the PMFCQ condition is essential for the validity of both normal cone representations in \eqref{4.1} and \eqref{4.2}; moreover, this condition cannot be replaced by its weaker EMFCQ version.

\begin{Example}\label{exa3} {\bf (violation of the normal cone representations with no PMFCQ).} {\rm Consider the infinite inequality system in
$\R^2$ given in Example~\ref{exa1}. It is shown therein that the EMFCQ holds at $\ox=(-1,0)$ while the PMFCQ does not. It is easy to check that in this case $\Hat N(\ox;\O)=N(\ox;\O)=\R_+\times\R_-$ while
\begin{equation*}
{\rm cl}\;{\rm cone}\big\{\nabla g_t(\ox)\big|\;t\in T_\ve(\ox)\big\}={\rm cl}\;{\rm cone}\big\{(1,0)\cup\{(t,0)|\;t\in(0,\ve)\big\}\subset\R_+\times\{0\}.
\end{equation*}
i.e., the inclusions ``$\subset$" in \eqref{4.1} and \eqref{4.2} are violated.}
\end{Example}

The next example shows that the perturbed active index set $T_\ve(\ox)$ cannot be replaced by its unperturbed counterpart $T(\ox)$ in the normal cone representations \eqref{4.1} and \eqref{4.2}.

\begin{Example}\label{ex3} {\bf (perturbation of the active index set is essential for the normal cone representations).}
{\rm Let us reconsider the nonlinear infinite system in problem (\ref{3.7}):
\begin{eqnarray*}
\left\{\begin{array}{ll}
g_1(x)=x_1+1\le 0,\\
g_n(x)=\disp\frac{1}{3n}x^3_1-x_2\le 0,\;\;n\in \IN\setminus\{1\},
\end{array}\right.
\end{eqnarray*}
where $x=(x_1,x_2)\in\R^2$ and $T:=\IN$. It is easy to check this inequality system satisfies our standing assumptions and that
the functions $g_t$  are uniformly strictly differentiable at $\ox=(-1,0)$. Observe further that $\O=\{(x_1,x_2)\in \R^2|\;x_1\le-1,
\;x_2\ge 0\}$ and hence $N(\ox;\O)=\R_+\times\R_-$. As shown above, both PMFCQ and EMFCQ conditions hold at $\ox$. However, we have
$T(\ox)=\{1\}$ and
$$
N(\ox;\O)\not=\cone\big\{\nabla g_t(\ox)\big|\;t\in T(\ox)\big\}=\cone\big\{\nabla g_1(\ox)\big\}=\cone\{(1,0)\}=\R_+\times \{0\},
$$
which shows the violation of the unperturbed counterparts of \eqref{4.1} and \eqref{4.2}. Observe that
\begin{eqnarray*}
\cone\big\{\nabla g_t(\ox)\big|\;t\in T_\ve(\ox)\big\}&=&\cone\Big\{(1,0)\cup\big\{\big(\frac{1}{n},-1\big)\big|\; n\in\IN\setminus\{1\},\;n
\ge\frac{1}{\ve}\Big\}\\
&=&\big\{(x_1,x_2)\in\R^2\big|\;x_1\ge 0,\;x_2<0\big\},
\end{eqnarray*}
which is not a closed subset. On the other hand, we have
$$
N(\ox;\O)=\bigcap_{\ve>0}{\rm cl}\;\cone\big\{\nabla g_t(\ox)\big|\;t\in T_\ve(\ox)\big\},
$$
which illustrates the validity of the normal cone representations in Theorem~\ref{thm41}.}
\end{Example}

Now we derive several consequences of Theorem~\ref{thm41}, which are of their independent interest. The first one concerns the case when the $\{\nabla g_t(\ox)|\;t\in T\}$ may not be bounded in $X^*$ as in our standing assumptions. It follows that the latter case can be reduced to the basic case of Theorem~\ref{thm41} with some modifications.

\begin{Corollary}\label{mc} {\bf (normal cone representation for infinite systems with unbounded gradients).} Considering the constraint system \eqref{3.1}, assume the following:

{\bf (a)} The functions {$g_t$, $t\in T$, are Fr\'echet differentiable at the point $\ox$ with $\|\nabla g_t(\ox)\|>0$ for all $t\in T$ and the mapping $h$ is strictly differentiable at $\ox$.

{\bf (b)} We have that $\disp\lim_{\eta\dn0}\tilde r(\eta)=0$, where $\tilde r(\eta)$ is defined by
\begin{eqnarray}\label{r1}
\tilde r(\eta):=\sup_{t\in T}\sup_\substack{x,x^\prime\in\B_{\eta}(\ox)}{x\ne x'}\frac{|g_t(x)-g_t(x^\prime)-\la\nabla g_t(\ox),x-x'\ra|}{\|\nabla g_t(\ox)\|\cdot\|x-x'\|}\;\mbox{ for all }\;\eta>0.
\end{eqnarray}

{\bf (c)} The operator $\nabla h(\ox)\colon X\to Y$ is surjective and for some $\ve>0$ there are $\tilde x\in X$ and $\sigma>0$ such that $\nabla h(\ox)\tilde x=0$ and that
\begin{eqnarray}\label{4.11}
\la\nabla g_t(\ox),\tilde x+x\ra\le 0\;\mbox{ whenever }\;\|x\|\le\sigma
\end{eqnarray}
for each $t\in\Tilde T_\ve(\ox):=\{t\in T|\ g_t(\ox)\ge-\ve\|\nabla g_t(\ox)\|\}$. Then the limiting normal cone to $\O$ at $\ox$ is computed by formula \eqref{4.2}.}
\end{Corollary}
{\bf Proof.} Define $\tilde g_t(x):=g_t(x)\|\nabla g_t(\ox)\|^{-1}$ for all $x\in X$ and $t\in T$ and observe that the feasible set $\O$ from \eqref{1.2} admits the representation
$$
\O=\big\{x\in X\big|\;\tilde g_t(x)\le 0,\;h(x)=0\big\}.
$$
Replacing $g_t$ by $\tilde g_t$ in Theorem~\ref{thm41}, we have that the functions $\{\tilde g_t\}$ and $h$ satisfy the standing assumptions (SA) as well as condition (\ref{3.4}) with the function \eqref{r1} instead of $r(\eta)$. Furthermore, it follows from (\ref{4.11}) that for some $\ve>0$ there are $\tilde x\in X$ and $\sigma>0$ satisfying $\nabla h(\ox)\tilde x=0$ and such that
\begin{eqnarray*}
\la\nabla\tilde g_t(\ox),\tilde x\ra\le-\sup_{x\in\B_\sigma(\ox)}\la\nabla\tilde g_t(\ox),x\ra=-\sigma\|\nabla\tilde g_t(\ox)\|\;\mbox{ whenever }\; t\in \Tilde T_\ve(\ox),
\end{eqnarray*}
which turns into $\la\nabla\tilde g_t(\ox),\tilde x\ra\le-\sigma$ for all $t\in\Tilde T_\ve(\ox)=\{t\in T|\;\tilde g_t(\ox)\ge-\ve\}$. Hence the PMFCQ condition holds for the functions $\tilde g_t$ and $h$ at $\ox$. It follows from Theorem~\ref{thm41} that
\begin{eqnarray*}
\begin{array}{ll}
N(\ox;\O)&=\disp\bigcap_{\ve>0}\cl\cone\big\{\nabla \tilde g_t(\ox)\big|\; t\in\Tilde T_\ve(\ox)\big\}+\nabla h(\ox)^*(Y^*)\\
&=\disp\bigcap_{\ve>0}\cl\cone\big\{\nabla g_t(\ox)\,\|\nabla g_t(\ox)\|^{-1}\big|\;t\in \Tilde T_\ve(\ox)\big\}+\nabla h(\ox)^*(Y^*)\\
&=\disp\bigcap_{\ve>0}\cl\cone\big\{\nabla g_t(\ox)\big|\;t\in\Tilde T_\ve(\ox)\big\}+\nabla h(\ox)^*(Y^*),
\end{array}
\end{eqnarray*}
which gives \eqref{4.2} and completes the proof of the corollary. $\h$\vspace*{0.1in}

Now we compare the result of Corollary~\ref{mc} with the recent one obtained in
\cite[Theorem~3.1 and Corollary~4.1]{S} for inequality constraint systems, i.e., with $h=0$ in \eqref{3.1}. The latter result is given by the inclusion form
\begin{eqnarray*}
N(\ox;\O)\subset\disp\bigcap_{\ve>0}\cl\cone\big\{\nabla g_t(\ox)\big|\;t\in T_\ve(\ox)\big\}
\end{eqnarray*}
in the case of $\|\nabla g_t(\ox)\|=1$ for all $t\in T$ under the Fr\'echet differentiability of $g_t$ {\em around} $\ox$ (in (as) we need it merely {\em at} $\ox$) and the replacement of (b) of Corollary~\ref{mc} by the following {\em equicontinuity} requirement on $g_t$ at $\ox$:
for each $\gamma>0$ there is $\eta>0$ such that
\begin{eqnarray}\label{4.12}
\|\nabla g_t(x)-\nabla g_t(\ox)\|\le\gamma\;\mbox{ for all }\;x\in\B_\eta(\ox),\;t\in T.
\end{eqnarray}
Let us check that the latter assumption together with the Fr\'echet differentiability of $g_t$ around $\ox$ imply (b) in Corollary~\ref{mc}.
Indeed, suppose that (\ref{4.12}) holds and then pick any $x,x^\prime\in\B_\eta(\ox)$. Employing the classical Mean Value Theorem, find $\hat x\in [x,x^\prime]\subset\B_\eta(\ox)$ such that $g_t(x)-g_t(x^\prime)=\la\nabla g_t(\hat x),x-x^\prime\ra$. This gives
\begin{eqnarray*}
\disp\frac{|g_t(x)-g_t(x')-\la\nabla g_t(\ox),x-x'\ra|}{\|\nabla g_t(\ox)\|\cdot\|x-x'\|}&=&\disp\frac{|\la\nabla g_t(\hat x),x-x^\prime\ra-\la\nabla g_t(\ox),x-x'\ra|}{\|x-x'\|}\\
&\le&\disp\frac{|\la\nabla g_t(\hat x)-\nabla g_t(\ox),x-x^\prime\ra|}{\|x-x'\|}\\
&\le&\disp \|\nabla g_t(\hat x)-\nabla g_t(\ox)\|\le\gamma
\end{eqnarray*}
and yields $\disp\lim_{\eta\dn 0}\tilde r(\eta)\le\gamma$ for all $\gamma>0$, which ensures the validity of (b) in Corollary~\ref{mc}.\vspace*{0.1in}

The next consequence of Theorem~\ref{thm41} concerns problems of semi-infinite programming and presents sufficient conditions for the fulfillment of simplified representations of the normal cones to feasible constraints with no closure operations in \eqref{4.1} and \eqref{4.2} and with the replacement of the perturbed index set $T_\ve(\ox)$ by that of active constraints $T(\ox)$.

\begin{Corollary}\label{si} {\bf (normal cones for  semi-infinite constraints).}  Let $X$ and $Y$ be finite-dimensional spaces with $\dim Y<\dim X$. Assume that $T$ is a compact metric space, that the function $t\in T\mapsto g_t(\ox)$ is u.s.c., and the mapping $t\in T\mapsto\nabla g_t(\ox)$ is continuous. Suppose further that system {\rm(\ref{3.1})} satisfies the PMFCQ at $\ox$. Then we have
\begin{eqnarray}\label{4.13}
\Tilde N(\ox;\O)=\disp\cone\big\{\nabla g_t(\ox)\big|\;t\in T(\ox)\big\}+\nabla h(\ox)^*(Y^*),
\end{eqnarray}
where $\Tilde N(\ox;\O)=\Hat N(\ox;\O)$ when the functions $g_t$ are uniformly Fr\'echet differentiable at $\ox$ and $\Tilde N(\ox;\O)= N(\ox;\O)$ when $g_t$ are uniformly strictly differentiable at $\ox$.\\[1ex]
In particular, if we assume in addition that both $t\in T\mapsto g_t(\ox)$ and $(x,t)\in X\times T\mapsto\nabla g_t(x)$ are continuous, then we also have {\rm(\ref{4.13})} for $\Tilde N(\ox;\O)= N(\ox;\O)$  provided that merely the EMFCQ condition holds at $\ox$.
\end{Corollary}
{\bf Proof.} Let $X=\R^d$ for some $d\in\IN$. It follows from Proposition~\ref{p3} that $g_t$, $t\in T$, and $h$ satisfy our standing assumptions $(SA)$. Since system (\ref{3.1}) satisfies the PMFCQ at $\ox$, there are $\tilde\ve>0$, $\delta>0$, and $\tilde x\in X$ such that $\la\nabla g_t(\ox),\tilde x\ra<-\delta$ for all $t\in T_\ve(\ox)$ and $\ve\in (0,\tilde \ve)$. Observe that the perturbed active index set $T_\ve(\ox)$ is compact in $T$ for all $\ve>0$ due to the u.s.c. assumption on $t\in T\mapsto g_t(\ox)$. It follows from the continuity of $t\in T\mapsto\nabla g_t (\ox)$ that $\{\nabla g_t(\ox)|\; t\in T_\ve(\ox)\}$ is a compact subset of $\R^d$.

We now claim that $0\notin\co\{\nabla g_t(\ox)|\;t\in T_\ve(\ox)\}$. Indeed, it follows for any $\lm\in\Tilde \R_+^{T_\ve(\ox)}$ with $\sum_{t\in T_\ve(\ox)}\lm_t=1$ that
\begin{eqnarray*}
\sum_{t\in T_\ve(\ox)}\lm_t\la\nabla g_t(\ox),\tilde x\ra\le-\sum_{t\in T_\ve(\ox)}\lm_t\delta=-\delta<0,
\end{eqnarray*}
which yields that $0\ne\sum_{t\in T_\ve(\ox)}\lm_t\nabla g_t(\ox)$, i.e., $0\notin\co\{\nabla g_t(\ox)|\;t\in T_\ve(\ox)\}$.

Hence it follows from \cite[Proposition 1.4.7]{HL} that the conic hull cone$\{\nabla g_t(\ox)|\;t\in T_\ve(\ox)\}$ is closed in $\R^d$. Combining this with Theorem~\ref{thm41}, it suffices to show that
\begin{eqnarray}\label{4.14}
\bigcap_{\ve>0}\cone\big\{\nabla g_t(\ox)\big|\;t\in T_\ve(\ox)\big\}=\cone\big\{\nabla g_t(\ox)\big|\;t\in T(\ox)\big\}.
\end{eqnarray}
Observe that the inclusion ``$\supset$'' in (\ref{4.14}) is obvious due to $T(\ox)\subset T_\ve(\ox)$ as $\ve>0$. To justify the converse inclusion, pick an arbitrary element $x^*$ from the set on the left-hand side of (\ref{4.14}). By the classical Carath\'eodory theorem, for all large $n\in\IN$ we find $\lm_n\in \R^{d+1}_+$ and
$$
\nabla g_{t_{n_1}}(\ox),\ldots,\nabla g_{t_{n_{d+1}}}(\ox)\in\big\{\nabla g_t(\ox)\big|\;t\in T_{\frac{1}{n}}(\ox)\big\}\subset\R^d
$$
satisfying the relationship
\begin{eqnarray}\label{4.15}
x^*=\sum_{k=1}^{d+1}\lm_{n_k}\nabla g_{t_{n_k}}(\ox),
\end{eqnarray}
which implies in turn that
\begin{eqnarray*}
\la x^*,\tilde x\ra=\sum_{k=1}^{d+1}\lm_{n_k}\la\nabla g_{t_{n_k}}(\ox),\tilde x\ra\le-\sum_{k=1}^{d+1}\lm_{n_k}\delta.
\end{eqnarray*}
Hence the sequence $\{\lm_n\}$ is bounded in $\R^{d+1}$, and so is
$$
\big\{\lm_n\times(\nabla g_{t_{n_1}}(\ox),\ldots,\nabla g_{t_{n_{d+1}}})\big\}\subset\R^{d+1}\times\R^{d(d+1)}.
$$
By the classical Bolzano-Weierstrass theorem and the compactness of $T$, we assume without loss of generality that the sequence $\{t_{n_k}\}$ converges to some $\bar t_k\in T$ for each $1\le k\le d+1$ and that $\{\lm_n\}$ converges to some $\bar\lm\in \R^{d+1}$ as $n\to \infty$. Note that $0\ge g_{t_{n_k}}(\ox)\geq -\frac{1}{n}$ for all $n\in\IN$ sufficiently large, which gives us
\begin{equation*}
0\ge g_{\bar t_k}(\ox)\ge\disp \limsup_{n\to\infty}g_{t_{n_k}}(\ox)\ge\limsup_{n\to\infty}-\frac{1}{n}=0
\end{equation*}
for all $1\le k\le d+1$. Combining the latter with (\ref{4.15}) ensures that
\begin{eqnarray*}
x^*=\sum_{k=1}^{d+1}\bar\lm_{k}\nabla g_{\bar t_k}(\ox)\in\cone\big\{\nabla g_t(\ox)\big|\;t\in T(\ox)\big\},
\end{eqnarray*}
which yields the inclusion ``$\subset$'' in (\ref{4.14}). Thus we arrive at formula $(\ref{4.13})$.

The second part of the corollary follows from the first part, Proposition~\ref{p3}, and Proposition~\ref{p4}. This completes the proof of the claimed result.  $\h$\vspace*{0.1in}

The results obtained in Corollary~\ref{si} can be compared with \cite[Theorem 3.4]{CHY}, where ``$\subset$'' in (\ref{4.13}) was obtained for $h=0$ under the following conditions: $T$ is scattered compact (meaning that every subset $S\subset T$ has an isolated point), $g_t$ are Fr\'echet differentiable for all $t\in T$, the mappings $(x,t)\in X\times T\mapsto g_t(x)$ and $(x,t)\in X\times T\mapsto\nabla g_t(x)$ are continuous, and the EMFCQ condition holds at $\ox$. We can see that these assumptions are significantly stronger than those Corollary~\ref{si}. Note, in particular, that the scattering compactness requirement on the index set $T$ is not different in applications from $T$ being finite.\vspace*{0.1in}

The next question we address in this section is about the possibility to obtain normal cone representations of the ``unperturbed" type as in Corollary~\ref{si} while in infinite programming settings with no finite dimensionality, compactness, and continuity assumptions made above. The following theorem shows that this can be done when the PMFCQ is accompanied by the NFMCQ condition of Definition~\ref{cd3}.

\begin{Theorem}\label{thm42} {\bf (unperturbed representations of normal cones for infinite constraint systems).} Let the functions $g_t$, $t\in T$, be uniformly Fr\'echet differentiable at $\ox$, and let that system {\rm(\ref{3.1})} satisfy the PMFCQ and NFMCQ conditions at $\ox$.
Then
\begin{eqnarray}\label{4.16}
\Hat N(\ox;\O)= \disp\cone\big\{\nabla g_t(\ox)\big|\;t\in T(\ox)\big\}+\nabla h(\ox)^*(Y^*).
\end{eqnarray}
If in addition the functions $g_t$, $t\in T$, are uniformly strictly differentiable at $\ox$, then
\begin{eqnarray}\label{4.17}
N(\ox;\O)=\disp\cone\big\{\nabla g_t(\ox)\big|\;t\in T(\ox)\big\}+\nabla h(\ox)^*(Y^*).
\end{eqnarray}
\end{Theorem}
{\bf Proof.} First we claim that the set $\disp\bigcap_{\ve>0}\cl\cone\{\nabla g_t(\ox)|\; t\in T_\ve(\ox)\}$ belongs to the set
\begin{eqnarray}\label{4.18}
\Big\{x^*\in X^*\Big|\;(x^*,\la x^*,\ox\ra)\in\cl\cone\big\{(\nabla g_t(\ox),\la\nabla g_t(\ox),\ox\ra-g_t(\ox))\big|\; t\in T\big\}\Big\}.
\end{eqnarray}
Indeed, it follows from the PMFCQ for (\ref{3.1}) at $\ox$ that $\nabla h(\ox)$ is surjective and  there are $\tilde \ve>0$, $\delta>0$, and $\tilde x\in X$ such that $\nabla h(\ox)\tilde x=0$ and that $\la\nabla g_t(\ox),\tilde x\ra<-\delta$ for all $\ve\le\tilde\ve$ and $t\in T_\ve(\ox)$. To justify the claimed inclusion to (\ref{4.18}), pick an arbitrary element $x^*\in\disp\bigcap_{\ve>0}\cl\cone\{\nabla g_t(\ox)|\;t\in T_\ve(\ox)\}$ and  for any $\ve\in (0,\tilde\ve)$ find a net $(\lm_\nu)_{\nu\in\N}\subset\tR$ with
\begin{eqnarray}\label{4.19}
x^*=w^*-\lim_\nu\disp\sum_{t\in T_\ve(\ox)}\lm_{t\nu}\nabla g_t(\ox).
\end{eqnarray}
This implies the relationships
\begin{eqnarray}\label{4.20}
\la x^*,\tilde x\ra=\lim_\nu\disp\sum_{t\in T_\ve(\ox)}\lm_{t\nu}\la \nabla g_t(\ox),\tilde x\ra\le-\delta\limsup_\nu\disp\sum_{t\in T_\ve(\ox)}\lm_{t\nu}\;\mbox{ and}
\end{eqnarray}
\begin{eqnarray*}
\begin{array}{ll}
\la x^*,\ox\ra=\disp\lim_\nu\disp\sum_{t\in T_\ve(\ox)}\lm_{t\nu}\la\nabla g_t(\ox),\ox\ra=\lim_\nu\disp\sum_{t\in T_\ve(\ox)}\lm_{t\nu}(\la\nabla g_t(\ox),\ox\ra-g_t(\ox)+g_t(\ox)).
\end{array}
\end{eqnarray*}
The later equality together with (\ref{4.20}) give us that
\begin{eqnarray*}
0\ge\la x^*,\ox\ra-\limsup_\nu\sum_{t\in T_\ve(\ox)}\lm_{t\nu}(\la\nabla g_t(\ox),\ox\ra-g_t(\ox))\ge\liminf_\nu\sum_{t\in T_\ve(\ox)}\lm_{t\nu} g_t(\ox)\ge\disp\frac{\ve}{\delta}\la x^*,\tilde x\ra.
\end{eqnarray*}
By passing to a subnet and combining this with (\ref{4.19}), we get
\begin{equation*}
(x^*,\la x^*,\ox\ra)\in\cl\cone\big\{(\nabla g_t(\ox),\la\nabla g_t(\ox),\ox\ra-g_t(\ox))\big|\;t\in T\big\}+\{0\}\times\disp[\frac{\ve}{\delta}\la x^*,\tilde x\ra,0]
\end{equation*}
for all $\ve\in(0,\tilde \ve)$, which implies that $x^*$ belongs to the set in (\ref{4.18}) by taking $\ve\dn 0$.

Involving further the NFMCQ condition, we claim the equality
\begin{eqnarray}\label{4.21}
\disp\bigcap_{\ve>0}\cl\cone\big\{\nabla g_t(\ox)\big|\;t\in T_\ve(\ox)\big\}=\cone\big\{\nabla g_t(\ox)\big|\;t\in T(\ox)\big\}.
\end{eqnarray}
The inclusion ``$\supset$'' in \eqref{4.21} is obvious since $T(\ox)\subset T_\ve(\ox)$ for all $\ve>0$. To justify the converse inclusion, pick any $x^*$ belonging to the left-hand side of (\ref{4.21}). By the NFMCQ condition, it follows from (\ref{4.18}) that there is $\lm\in\tR$  such that
\begin{eqnarray}\label{4.22}
(x^*,\la x^*,\ox\ra)=\sum_{t\in T}\lm_t\big(\nabla g_t(\ox),\la\nabla g_t(\ox),\ox\ra-g_t(\ox)\big),
\end{eqnarray}
which readily yields the equalities
\begin{eqnarray*}
0=\sum_{t\in T}\lm_t\la\nabla g_t(\ox),\ox\ra-\sum_{t\in T}\lm_t\big(\la\nabla g_t(\ox),\ox\ra-g_t(\ox)\big)=\sum_{t\in T}\lm_t g_t(\ox).
\end{eqnarray*}
Since $g_t(\ox)\le 0$, we get $\lm_tg_t(\ox)=0$ for all $t\in T$. Combining this with (\ref{4.22}) gives us
\begin{eqnarray*}
x^*\in\cone\big\{\nabla g_t(\ox)\big|\;t\in T(\ox)\big\},
\end{eqnarray*}
which implies the inclusion ``$\subset$'' in (\ref{4.21}). To complete the proof of the theorem, we combine the obtained equality \eqref{4.21} with finally Theorem~\ref{thm41}. $\h$\vspace*{0.1in}

Observe from Proposition~\ref{fm} that formula (\ref{4.16}) holds  under our standing assumptions (SA) and the MFCQ condition at $\ox$ when $T$ is a finite index set. Furthermore, the formula for the limiting normal cones (\ref{4.17}) is  also satisfied if all the functions $g_t$ are strictly differentiable at $\ox$. It follows from  Proposition~\ref{fm} that Corollary~\ref{si} can be derived from a semi-infinite version of Theorem~\ref{thm42} in addition to the assumptions of this corollary we suppose that the function $t\in T\mapsto g_t(\ox)$ is continuous in $T$.\vspace*{0.05in}

The next example shows that the PMFCQ condition cannot be replaced by the EMFCQ one in Theorem~\ref{thm42} to ensure the unperturbed normal cone representations \eqref{4.16} and \eqref{4.17} in the presence of the NFMCQ.

\begin{Example}\label{ex4}{\bf (EMFCQ combined with NFMCQ does not ensure the unperturbed normal cone representations).}
{\rm We revisit the semi-infinite inequality constraint system in Example~\ref{ex1}. It is shown there that this system satisfied the EMFCQ but not PMFCQ at $\ox=(-1,0)$. It is easy to check that the set
\begin{eqnarray*}
\cone\big\{\big(\nabla g_t(\ox),\la\nabla g_t(\ox),\ox\ra-g_t(\ox)\big)\big|\;t\in T\big\}&=&\cone\big((1,0,-1)\cup\{(t,0,0)\big|t\in (0,1]\}\big)\\
&=&\big\{x\in\R^3\big|\;x_1+x_3\ge 0,\;x_1\ge 0\ge x_3,\;x_2=0\big\}
\end{eqnarray*}
is closed in $\R^3$, i.e., the NFMCQ condition holds at $\ox$. Observe however that both representations (\ref{4.16}) and (\ref{4.17}) are not satisfied for this system since we have
\begin{eqnarray*}
\Hat N(\ox;\O)=N(\ox;\O)\not=\cone\big\{\nabla g_t(\ox)\big|\;t\in T(\ox)\big\}=\cone\{(1,0)\}=\R_+\times\{0\}.
\end{eqnarray*}}
\end{Example}

Now we present a consequence of Theorem~\ref{thm42} with the corresponding discussions.

\begin{Corollary}\label{cv} {\bf (normal cone for infinite convex systems).}  Assume that all the functions $g_t$, $t\in T$, in \eqref{3.1} are convex and uniformly Fr\'echet differentiable and that $h=A$ is a surjective continuous linear operator. Suppose further that system {\rm(\ref{3.1})} satisfies the PMFCQ $($equivalently the SSC$)$ at $\ox\in\O$. Then the normal cone to $\O$ at $\ox$ in sense of convex analysis is computed by
\begin{eqnarray*}
N(\ox;\O)=\disp\bigcap_{\ve>0}\cl\cone\big\{\nabla g_t(\ox)\big|\;t\in T_\ve(\ox)\big\}+A^*(Y^*).
\end{eqnarray*}
If in addition the NFMCQ holds at $\ox$, then we have
\begin{eqnarray}\label{cnm}
N(\ox;\O)=\disp\cone\big\{\nabla g_t(\ox)\big|\;t\in T(\ox)\big\}+A^*(Y^*).
\end{eqnarray}
\end{Corollary}
{\bf Proof.} It follows directly from Proposition~\ref{ss} and Theorem~\ref{thm42}.  $\h$\vspace*{0.1in}

For $h=0$ in \eqref{3.1} the equality in (\ref{cnm}) can be deduced from \cite[Corollary~3.6]{DMN2} under another {\em Farkas-Minkowski Constraint Qualification} (FMCQ) defined as follows:\vspace*{0.05in}

{\bf (FMCQ)} {\em The conic hull {\rm cone}$\{\epi g_t^*|\; t\in T\}$ is weak$^*$ closed in $X^*\times\R$ under the additional assumption that the functions $g_t$ are l.s.c., where
\begin{equation*}
\ph^*(x^*):=\sup\big\{\la x^*,x\ra-\ph(x)\big|\;x\in X\big\},\quad x^*\in X^*,
\end{equation*}
stands for the Fenchel conjugate of a convex function.}\vspace*{0.1in}

It is worth noting that the above FMCQ condition is a {\em global} property, and hence formula (\ref{cnm}) holds at every $\ox\in\O$. By the contrary, our new NFMCQ condition (\ref{fm}) is constructed at a fixed point $\ox\in\O$. The next example shows that the combination of the PMFCQ (or the SSC) and the NFMCQ conditions for infinite convex inequality systems is not stronger than the FMCQ one.

\begin{Example}\label{ex5} {\bf (PMFCQ combined with NFMCQ does not imply FMCQ for convex inequality systems).} {\rm Define a function $g_t:\R^2\to\R$ by $g_t(x_1,x_2):=tx_1^2-x_2$ for all $(x_1,x_2)\in\R^2$ and $t\in T:=(0,1)$, and let $\ox=(0,0)\in\R^2$. It is easy to see that all the functions $g_t$, $t\in T$, are convex and differentiable and that the standing assumptions are satisfied. For each $t\in T$ we have
\begin{displaymath}
g_t^*(a_1,a_2)=\sup_{(x_1,x_2)\in\R^2}\left\{a_1x_1+a_2x_2-tx_1^2+x_2\right\}=\left\{
\begin{array}{ll}
\disp\frac{a_1^2}{4t} &\mbox{if $a_2=-1$},\\
\infty &\mbox{otherwise}.
\end{array}
\right.
\end{displaymath}
This implies that $\epi g_t^*=\{(a,-1,\frac{a^2}{4t}+r)|\;a\in \R,\;r\ge 0\}$, which yields in turn that
$$
C:=\cone\big\{\epi g_t^*\big|\;t\in T\big\}=\cone\Big\{(a,-1,\frac{a^2}{4}+r)\Big|\;a\in \R,\;r\ge 0\Big\}.
$$
The latter set is not closed in $\R^3$ since $\{0\}\times\{0\}\times\R_+\not\subset C$ while $\{0\}\times\{0\}\times \R_+\subset{\rm cl}C$. Moreover,  we see that $\nabla g_t(\ox)=(0,-1)$ for all $t\in T$, and then the PMFCQ  is satisfied. Furthermore, it follows that the set
\begin{equation*}
\cone\big\{(\nabla g_t(\ox),\la\nabla g_t(\ox),\ox\ra-g_t(\ox))\big|\;t\in T\big\}=\cone\{(0,-1, 0)\}=\{0\}\times\R_-\times\{0\}
\end{equation*}
is closed in $\R^3$. Hence the PMFCQ and NFMCQ conditions hold but the FMCQ does not.}
\end{Example}

Finally in this section, we give specifications of obtained normal cone representations in the case linear infinite systems.

\begin{Proposition}\label{clm} {\bf (normal cone representations for linear infinite constraint systems).} Consider the constraint system \eqref{3.1} with $g_t(x)=\la a^*_t,x\ra-b_t$ for all $t\in T$, and let $h=A:X\to Y$. Assume that $A$ is a surjective continuous linear operator and that the coefficient set $\{a^*_t|\;t\in T\}$ is bounded in $X^*$. If the SSC condition holds at $\ox$, then
\begin{eqnarray*}
N(\ox;\O)=\disp\bigcap_{\ve>0}\cl\cone\big\{a^*_t\big|\;t\in T_\ve(\ox)\big\}+A^*(Y^*)
\end{eqnarray*}
for the feasible set $\O:=\{x\in X|\;Ax=0,\;\la a^*_t,x\ra-b_t\le 0,\;t\in T\}$. On the other hand, assuming the weak$^*$ closedness of $\cone\{(a^*_t, b_t)|\;t\in T\}$ in $X^*\times\R$ and that $h=0$ gives us
\begin{eqnarray*}
N(\ox;\O)=\disp\cone\big\{a^*_t\big|\;t\in T(\ox)\big\}.
\end{eqnarray*}
\end{Proposition}
{\bf Proof.} The first statement is a specification of Corollary~\ref{cv}. The second one follows from the proofs given in \cite[Proposition~3.1]{CLMP1} and \cite[Theorem~3.2]{CLMP2} by using the classical Farkas Lemma for linear infinite systems. $\h$

\section{Optimality Conditions in Nonlinear Infinite Programming}
\setcounter{equation}{0}

In this section we employ general principles in optimization and the calculus results on computing the normal cones to the infinite constraint sets in Section~4 to deriving necessary optimality conditions for problems of infinite and semi-infinite programming. We confine ourselves to optimality conditions of the ``lower" subdifferential type conventional in minimization. Condition of the other (``upper" or superdifferential) type can be derived from the calculus results of Section~4 using an approach developed in \cite[Chapter~5]{M1}; see also the recent paper \cite{CLMP2} for the implementation of the latter approach in the case of semi-infinite and infinite programs with linear constraints.\vspace*{0.05in}

Our first theorem in this section concerns infinite programs of type \eqref{1.1} in arbitrary Banach spaces involving Fr\'echet differentiable cost functions.

\begin{Theorem}\label{thm43} {\bf (necessary optimality conditions for differentiable infinite programs in general Banach spaces).} Let $\ox$ be a local minimizer of the infinite program {\rm(\ref{1.1})} under the PMFCQ condition imposed on the constraints at $\ox$. Suppose further that the inequality constraint functions $g_t$, $t\in T$, are uniformly Fr\'echet differentiable at $\ox$ and the cost function $f$ is Fr\'echet differentiable at this point. Then we have the inclusion
\begin{eqnarray}\label{5.1}
0\in\nabla f(\ox)+\disp\bigcap_{\ve>0}\cl\cone\big\{\nabla g_t(\ox)\big|\;t\in T_\ve(\ox)\big\}+\nabla h(\ox)^*(Y^*).
\end{eqnarray}
If in addition the NFMCQ holds at $\ox$, then there exist multipliers $\lm\in\tR$ and $y^*\in Y^*$ satisfying the differential KKT condition
\begin{eqnarray}\label{5.2}
0=\nabla f(\ox)+\disp\sum_{t\in T(\ox)}\lm_t\nabla g_t(\ox)+\nabla h(\ox)^*y^*.
\end{eqnarray}
\end{Theorem}
{\bf Proof.} It is clear that $\ox$ is a local optimal solution to the following unconstrained optimization problem with the {\em infinite penalty}:
\begin{eqnarray}\label{5.3}
\mbox{minimize }\;f(x)+\delta(x;\O),
\end{eqnarray}
where $\O$ is the feasible constraint set \eqref{1.2}. Applying the generalized Fermat rule to the latter problem (see, e.g., \cite[Proposition~1.114]{M1}), we have
\begin{eqnarray}\label{5.4}
0\in\Hat\partial\big(f+\delta(\cdot;\O)\big)(\ox).
\end{eqnarray}
Since $f$ is Fr\'echet differentiable at $\ox$, it follows from the sum rule of \cite[Theorem~1.107 ]{M1} applied to \eqref{5.4} and from the first relationship in (\ref{2.4}) that
\begin{equation}\label{5.4a}
0\in\nabla f(\ox)+\Hat\partial\delta(\ox;\O)(\ox)=\nabla f(\ox)+\Hat N(\ox;\O).
\end{equation}
 Now using the Fr\'echet normal cone representation of Theorem~\ref{thm41} in \eqref{5.4a}, we arrive at (\ref{5.1}). The second part (\ref{5.2}) of this theorem readily follows from Theorem~\ref{thm42}. $\h$\vspace*{0.1in}

The next theorem establishes necessary conditions for local minimizers of infinite programs (\ref{1.1}) with general nonsmooth cost functions in the framework of Asplund spaces.

\begin{Theorem}\label{thm44} {\bf (necessary optimality conditions for nonconvex infinite programs defined on Asplund spaces, I).} Let $\ox$ be a local minimizer of problem {\rm(\ref{1.1})}, where the domain space $X$ is Asplund while the image space $Y$ is arbitrary Banach. Suppose that the constraint functions $g_t$, $t\in T$, are uniformly strictly differentiable at $\ox$, that the cost function $f$ is l.s.c.\ around $\ox$ and SNEC at this point, and that the qualification condition
\begin{eqnarray}\label{5.5}
\partial^\infty f(\ox)\cap\Big[-\disp\bigcap_{\ve>0}\cl\cone\big\{\nabla g_t(\ox)\big|\;t\in T_\ve(\ox)\big\}-\nabla h(\ox)^*(Y^*)\Big]=\{0\}
\end{eqnarray}
is fulfilled; the latter two assumptions are automatic when $f$ is locally Lipschitzian around $\ox$. If the PMFCQ condition holds at $\ox$, then
\begin{eqnarray}\label{5.6}
0\in\partial f(\ox)+\disp\bigcap_{\ve>0}\cl\cone\big\{\nabla g_t(\ox)\big|\;t\in T_\ve(\ox)\big\}+\nabla h(\ox)^*(Y^*).
\end{eqnarray}
If in addition we assume that the NFMCQ holds at $\ox$ and replace \eqref{5.5} by
\begin{eqnarray}\label{5.5a}
\partial^\infty f(\ox)\cap\Big[-\cone\big\{\nabla g_t(\ox)\big|\;t\in T(\ox)\big\}-\nabla h(\ox)^*(Y^*)\Big]=\{0\},
\end{eqnarray}
then there exist multipliers $\lm\in\tR$ and $y^*\in Y^*$ such that the following subdifferential KKT condition is satisfied:
\begin{eqnarray}\label{5.6a}
0\in\partial f(\ox)+\disp\sum_{t\in T(\ox)}\lm_t\nabla g_t(\ox)+\nabla h(\ox)^*y^*.
\end{eqnarray}
\end{Theorem}
{\bf Proof.} Observe first that the feasible set $\O$ is locally closed around $\ox$. Indeed, it follows from (\ref{3.4}) that there are $\gg>0$ and $\eta>0$ sufficiently small such that
\begin{eqnarray*}
\|h(x)-h(x^\prime)\|\le(\|\nabla h(\ox)\|+\gg)\|x-x^\prime\|\;\mbox{ and }\;\|g_t(x)-g_t(x^\prime)\|\le\sup_{\tau\in T}(\|\nabla g_\tau(\ox)\|+\gg)\|x-x^\prime\|
\end{eqnarray*}
for all $x,x^\prime\in\B_\eta(\ox)$ and $t\in T$. Picking any sequence $\{x_n\}\subset\O\cap\B_\eta(\ox)$ converging to some $x_0$ as $n\to\infty$, we have
\begin{eqnarray*}
\|h(x_0)\|\le(\|\nabla h(\ox)\|+\gg)\|x_n-x_0\|\;\mbox{ and }\;g_t(x_0)\le\sup_{\tau\in T}(\|\nabla g_\tau(\ox)\|+\gg)\|x_n-x_0\|+g_t(x_n)
\end{eqnarray*}
for each $t\in T$ and $n\in\IN$. By passing to the limit as $n\to\infty$, the latter yields that $h(x_0)=0$ and $g_t(x_0)\le 0$ for all $t\in T$, i.e., $x_0\in\O\cap\B_\eta(\ox)$, which justifies the local closedness of the feasible set $\O$ around $\ox$.

Employing now the generalized Fermat rule to the solution $\ox$ of \eqref{5.3} with the closed set $\O$ and using \cite[Theorem~3.36]{M1} on the sum rule for basic/limiting subgradients in Asplund spaces when $f$ is SNEC at $\ox$ yield that
\begin{equation}\label{sr}
0\in\partial\big(f+\delta(\cdot;\O)\big)(\ox)\subset\partial f(\ox)+\partial\delta(\ox;\O)=\partial f(\ox)+N(\ox;\O)
\end{equation}
provided that $\partial^\infty f(\ox)\cap\big(-N(\ox;\O)\big)=\{0\}$. We apply further to both latter conditions the limiting normal cone representation of Theorem~\ref{thm41}. This gives us the optimality condition \eqref{5.6} under the fulfillment of \eqref{5.5} and the PMFCQ at $\ox$. Applying finally Theorem~\ref{thm42} instead of Theorem~\ref{thm41} in the setting above, we arrive at the KKT condition \eqref{5.6a} under the assumed NFMCQ at $\ox$ and \eqref{5.5a}, which completes the proof of the theorem. $\h$\vspace*{0.1in}

An important ingredient in the proof of Theorem~\ref{thm44} is applying the subdifferential sum rule from \cite[Theorem~3.36 ]{M1} to the sum $f+\dd(\cdot;\O)$, which requires that either $f$ is SNEC at $\ox$ or $\O$ is SNC at this point. While the first possibility was used above, now we are going to explore the second alternative. The next proposition presents verifiable conditions ensuring the SNC property of the feasible set $\O$ at $\ox$.

\begin{Proposition}\label{asp} {\bf (SNC property of feasible sets in infinite programming).} Let $X$ be an Asplund space, and let $\dim Y<\infty$ in the framework of \eqref{1.1}. Assume that all the functions $g_t$, $t\in T$, are Fr\'echet differentiable around some $\ox\in\O$ and that the corresponding derivative family $\{\nabla g_t\}_{t\in T}$ is equicontinuous around this point, i.e., there exists $\ve>0$ such that for each $x\in\B_\ve(\ox)$ and each $\gg>0$ there is $0<\tilde\ve<\ve$ with the property
\begin{eqnarray}\label{5.7}
\|\nabla g_t(x^\prime)-\nabla g_t(x)\|\le\gg\;\mbox{ whenever }\;x^\prime\in\B_{\tilde\ve}(x)\cap \O\;\mbox{ and }\;t\in T.
\end{eqnarray}
Then the feasible set $\O$ in \eqref{1.2} is locally closed around $\ox$ and SNC at this point provided that the PMFCQ condition holds at $\ox$.
\end{Proposition}
{\bf Proof.} Consider first the set $\O_1:=\{x\in X|\;g_t(x)\le 0,\;t\in T\}$.  By using arguments similar to the proof of Theorem~\ref{thm44}, we justify the local closedness of $\O_1$ around $\ox$. Now let us prove that $\O_1$ is SNC at this point. To proceed, pick any sequence $(x_n,x^*_n)\in\O_1\times X^*$, $n\in\IN$, satisfying
\begin{eqnarray*}
x_n\st{\O_1}\to\ox,\;x^*_n\in\Hat N(x_n;\O_1)\;\mbox{ and }\;x^*_n\st{w^*}\to 0\;\mbox{ as }\;n\to\infty.
\end{eqnarray*}
Taking \eqref{5.7} into account,  we see that the functions $g_t$, $t\in T$ satisfy the standing assumptions (SA) at $x_n$ for all $n\in\IN$ sufficiently large. Moreover, the proof showing that assumption (\ref{3.4}) holds at $x_n$ follows from the discussions right after Corollary~\ref{mc}. Since the PMFCQ condition holds at $\ox$, there exist $\delta>0$, $\ve>0$, and $\tilde x\in X$ such that $\la\nabla g_t(\ox),\tilde x\ra\le-2\delta$ for all $t\in T_{2\ve}(\ox)$. Observe that $T_\ve(x_n)\subset T_{2\ve}(\ox)$ for all large $n\in\IN$. Indeed, whenever $t\in T_\ve(x_k)$ we have
\begin{eqnarray*}
0\ge g_t(\ox)&\ge& g_t( x_n)-\la \nabla g_t(\ox),x_n-\ox\ra-\|x_n-\ox\|s(\|x_k-\ox\|)\\
&\ge&-\ve-\sup_{\tau\in T}\|\nabla g_\tau(\ox)\|\cdot\|x_n-\ox\|-\|\tilde x_n-\ox\|s(\|x_n-\ox\|)\ge-2\ve
\end{eqnarray*}
for all large $n\in\IN$, where $s(\cdot)$ is defined in \eqref{3.3}. Further, it follows from (\ref{5.7}) that
\begin{eqnarray*}
\la\nabla g_t( x_n),\tilde x\ra\le\la\nabla g_t(\ox),\tilde x\ra+\|\nabla g_t( x_n)-\nabla g_t(\ox)\|\cdot\|\tilde x\|\le-2\delta+\|\nabla g_t( x_n)-\nabla g_t(\ox)\|\cdot\|\tilde x\|\le-\delta
\end{eqnarray*}
when $n\in\IN$ is sufficiently large. Hence we suppose without loss of generality that
\begin{eqnarray}\label{5.8}
T_\ve(x_n)\subset T_{2\ve}(\ox)\;\mbox{ and }\;\sup_{t\in T_\ve( x_n)}\la\nabla g_t(x_n),\tilde x\ra\le-\delta\;\mbox{ whenever }\;n\in\IN.
\end{eqnarray}
Applying now Theorem~\ref{thm41} in this setting, we have that for each $n\in\IN$ there exists a net $\{\lm_{n_\nu}\}_{\nu\in\N}\subset\Tilde \R_+^{T_\ve( x_n)}$ such that
\begin{eqnarray*}
x^*_n=w^*-\lim_\nu\sum_{t\in T_\ve(x_n)}\lm_{tn_\nu}\nabla g_t(x_n).
\end{eqnarray*}
Combining this with (\ref{5.8}) yields that
\begin{eqnarray*}
\la x^*_n,\tilde x\ra=\lim_\nu\sum_{t\in T_\ve(x_n)}\lm_{tn_\nu}\la\nabla g_t(x_n),\tilde x\ra\le-\delta\liminf_\nu\sum_{t\in T_\ve( x_n)}\lm_{tn_\nu}.
\end{eqnarray*}
Furthermore, for each $x\in X$ we get the relationships
\begin{eqnarray*}
\la x^*_n,x\ra&=&\liminf_\nu\sum_{t\in T_\ve( x_k)}\lm_{tn_\nu}\la\nabla g_t(x_n),x\ra\le\liminf_\nu\sum_{t\in T_\ve(x_n)}\lm_{tn_\nu}\sup_{\tau\in T}\|\nabla g_\tau(x_n)\|\cdot\|x\|\\
&\le&-\frac{\la x^*_n,\tilde x\ra}{\delta}\sup_{\tau\in T}\|\nabla g_\tau(x_n)\|\cdot\|x\|,
\end{eqnarray*}
which imply that $\|x^*_n\|\le-\frac{\la x^*_n,\tilde x\ra}{\delta}\sup_{\tau\in T}\|\nabla g_\tau(x_n)\|$ for all $n\in\IN$. Since $x^*_n\st{w^*}\to 0$, it follows from the latter that $\|x^*_n\|\to 0$ as $n\to\infty$ and thus the set $\O_1$ is SNC at $\ox$.

Consider now the set $\O_2:=\{x\in X|\;h(x)=0\}$, which is obviously closed around $\ox$. It follows from \cite[Theorem~1.22]{M1} and finite dimensionality of $Y$ that $\O_2$ is SNC at $\ox$. Moreover, we get from \cite[Theorem~1.17]{M1} that $N(\ox;\O_2)=\nabla h(\ox)^*(Y^*)$. Thus for any $x^*\in N(\ox;\O_1)\cap(-N(\ox;\O_2))$ there is $y^*\in Y^*$ such that $x^*+\nabla h(\ox)^*y^*=0$, and then
$$
\la x^*,\tilde x\ra=-\la\nabla h(\ox)^*y^*,\tilde x\ra=-\la y^*,\nabla h(\ox)\tilde x\ra=0.
$$
Since $x^*\in N(\ox;\O_1)$, we find by Theorem~\ref{thm41} such a net $\{\lm_\nu\}_{\nu\in\N}\in\tR$ that
\begin{eqnarray*}
x^*=w^*-\lim_\nu\sum_{t\in T_\ve(\ox)}\lm_{t\nu}\nabla g_t(\ox),
\end{eqnarray*}
which yields in turn that
\begin{eqnarray*}
0=\la x^*,\tilde x\ra=\lim_\nu\sum_{t\in T_\ve(\ox)}\lm_{t\nu}\la \nabla g_t(\ox),\tilde x\ra\le-2\delta\liminf_{\nu}\sum_{t\in T_\ve(\ox)} \lm_{t\nu}.
\end{eqnarray*}
This ensures the relationships
\begin{equation*}
\la x^*,x\ra=\liminf_\nu\sum_{t\in T_\ve(\ox)}\lm_{t\nu}\la\nabla g_t(\ox),x\ra\le\liminf_\nu\sum_{t\in T_\ve(\ox)}\lm_{t\nu}\sup_{\tau\in T}\|\la \nabla g_\tau(\ox)\|\|x\|=0
\end{equation*}
for all $x\in X$. Hence we have $x^*=0$, and so $N(\ox;\O_1)\cap(-N(\ox;\O_2))=\{0\}$. It finally follows from \cite[Corollary~3.81]{M1} that the intersection $\O=\O_1\cap\O_2$ is SNC at $\ox$, which thus completes the proof of the proposition. $\h$\vspace*{0.05in}

 Observe that the assumption $\dim Y<\infty $ is essential in Proposition~\ref{asp}. To illustrate this, consider a particular case of \eqref{1.1} when $T=\emp$. It follows from \cite[Theorem~1.22]{M1} that the inverse image $\O=h^{-1}(0)$ is SNC at $\ox\in\O$ if and only if the set $\{0\}$ is SNC at $0\in Y$. Since $N(0;\{0\})=Y^*$, the latter holds if and only if the weak$^*$ topology in $Y^*$ agrees with the norm topology in $Y^*$, which is only the case of $\dim Y<\infty$ by the classical Josefson-Nissenzweig theorem from theory of Banach spaces.\vspace*{0.05in}

 Now we are ready to derive an aforementioned alternative counterpart of Theorem~\ref{thm44}.

\begin{Theorem}\label{thm45} {\bf (necessary optimality conditions for nonconvex infinite programs defined on Asplund spaces, II).}
Let $\ox$ be a local minimizer of infinite program {\rm (\ref{3.1})} under the assumptions of Proposition~{\rm\ref{asp}}. Suppose also that $f$ is l.s.c.\ around $\ox$  and that the qualification condition \eqref{5.5} is satisfied. Then we have the optimality condition \eqref{5.6}. If in addition we assume that the NFMCQ holds at $\ox$ and replace \eqref{5.5} by \eqref{5.5a},
then there exist multipliers $\lm\in\tR$ and $y^*\in Y^*$ such that the subdifferential KKT condition \eqref{5.6a}.
\end{Theorem}
{\bf Proof.} It is similar to the proof of Theorem~\ref{thm44} with applying Proposition~\ref{asp} on the SNC and closedness property of $\O$ in the sum rule \eqref{sr} of \cite[Theorem 3.36]{M1}. $\h$\vspace*{0.05in}

The next result provides necessary and sufficient optimality conditions for convex problems of infinite programming in general Banach spaces.

\begin{Theorem}\label{conv1} {\bf (necessary and  optimality conditions for convex infinite programs).} Let both spaces $X$ and $Y$ be Banach. Assume that all the functions $g_t$, $t\in T$, are convex and uniformly Fr\'echet differentiable and that $h=A$ is a surjective  continuous linear operator. Suppose further that the cost function $f$ is convex and continuous at some point in $\O$. If the PMFCQ condition $($equivalently the SSC condition) holds at $\ox$, then $\ox$ is a global minimizer of problem \eqref{1.1} if and only if
\begin{eqnarray*}
0\in\partial f(\ox)+\disp\bigcap_{\ve>0}\cl\cone\big\{\nabla g_t(\ox)\big|\;t\in T_\ve(\ox)\big\}+A^*(Y^*).
\end{eqnarray*}
If in addition the NFMCQ condition holds, then $\ox$ is a global minimizer of problem \eqref{1.1} if and only if there exist $\lm\in \tR$ and $y^*\in Y^*$ such that
\begin{eqnarray}\label{5.9}
0\in\partial f(\ox)+\disp\sum_{t\in T(\ox)}\lm_t\nabla g_t(\ox)+A^*y^*.
\end{eqnarray}
\end{Theorem}
{\bf Proof.} Observe that $\ox$ is a global minimizer of problem (\ref{1.1}) if and only if it is a global minimizer of the convex unconstrained problem (\ref{5.3}), which is equivalent to the fact that
\begin{eqnarray*}
0\in\partial\big(f+\delta(\cdot;\O)\big)(\ox).
\end{eqnarray*}
Applying the convex subdifferential sum rule to the latter inclusion, we conclude that $\ox$ is a global minimizer of problem (\ref{1.1}) if and only if
\begin{eqnarray*}
0\in\partial f(\ox)+\partial\delta(\ox;\O)=\partial f(\ox)+N(\ox;\O).
\end{eqnarray*}
The rest of the proof follows from Corollary~\ref{cv}. $\h$\vspace*{0.1in}

Note that some versions of necessary optimality condition of the KKT type \eqref{5.9} were derived in \cite[Theorems~3.1 and 3.2]{CLMP2} for infinite problems with linear constraints but possibly nonconvex cost functions under the SSC and the linear counterpart of the FMCQ; see Example~\ref{ex5} and the corresponding discussions above.

Observe also that the results of Theorem~\ref{thm45} and Theorem~\ref{conv1} are formulated with no change in the case of semi-infinite programs, while in Theorem~\ref{thm43} we just drop the SNEC assumption on $f$, which holds automatically when $X$ is finite-dimensional.\vspace*{0.1in}

In conclusion we present a consequence of our results for the classical framework of semi-infinite programming while involving nonsmooth cost functions.

\begin{Corollary}\label{cor6} {\bf (necessary optimality conditions for semi-infinite programs with compact index sets).} Let $\ox$ be a local  minimizer of program {\rm(\ref{1.1})}, where both spaces $X$ and $Y$ are finite-dimensional with $\dim Y<\dim X$. Assume that the index set $T$ in \eqref{1.1} is a compact metric space, that the mappings $(x,t):X\times T\mapsto g_t(x)$ and $(x,t):X\times T\mapsto\nabla g_t(x)$ are continuous, and that the cost function $f$ is l.s.c.\ around $\ox$ with the fulfillment of \eqref{5.5a}. If in addition the EMFCQ holds at $\ox$, then there exist multipliers $\lm\in\tR$ and $y^*\in Y^*$ satisfying the subdifferential KKT condition \eqref{5.6a}.
\end{Corollary}
{\bf Proof.} By Proposition~\ref{sfm} we have that the NFMCQ condition holds at $\ox$ under the assumptions made. Then this corollary follows directly from Theorem~\ref{thm44}. $\h$\vspace*{0.1in}

When $f$ is smooth around $\ox$, assumption \eqref{5.5a} holds automatically while \eqref{5.6a} reduced to the differential KKT condition \eqref{5.2}. Then Corollary~\ref{cor6} reduces to a well-known result in semi-infinite programming that can be found, e.g., in \cite[Theorem~3.3]{HK} and \cite[Theorem~2]{LS}.

\end{document}